\documentclass[runningheads]{llncs}
\usepackage{amssymb,amsmath}
\usepackage{graphicx}
\usepackage{enumerate}
\usepackage{comment}
\usepackage{epstopdf}
\usepackage[bookmarksnumbered=true,hidelinks]{hyperref}
\usepackage{appendix}
\usepackage[numbers,sort&compress]{natbib}
\usepackage{etoolbox}
\usepackage{lipsum}% just to generate text for the example
\usepackage{todonotes}
\patchcmd{\thebibliography}{\chapter*}{\section*}{}{}
\usepackage{todonotes,bm}
\usepackage{url}
\usepackage{cleveref}
\usepackage[subtle, indent=normal, mathspacing=normal]{savetrees}
\usepackage{algorithm}
\usepackage{algorithmicx}
\usepackage{algpseudocode}
\algdef{SE}[DOWHILE]{Do}{doWhile}{\algorithmicdo}[1]{\algorithmicwhile\ #1}%

\def\Re{\mathbb{R}}

\def\hat{\widehat}

\def\N{{\mathcal N}}

\def\Re{{\mathbb R}}
\def\Z{{\mathcal Z}}

\newcommand{\exclude}[1]{}

\DeclareMathOperator{\Diag}{Diag}

\DeclareMathOperator{\rank}{rank}
\DeclareMathOperator{\tr}{tr}
\DeclareMathOperator{\argmax}{argmax}

\renewcommand{\S}{\mathcal{S}}
\newtheorem{claim1}{Claim}
\newtheorem{observation}{Observation}
\newcommand*{\qedA}{\hfill\ensuremath{\diamond}}

\allowdisplaybreaks

\begin{document}
	\mainmatter  % start of an individual contribution
	
	\title{On Sparse Canonical Correlation Analysis}
	
	\titlerunning{On  Sparse Canonical Correlation Analysis}
	
	
	\author{Yongchun Li, Santanu S. Dey, and Weijun Xie}

	\authorrunning{Y. Li, S. Dey, and W. Xie}
	
	\institute{ School of Industrial and Systems Engineering, Georgia Institute of Technology%\\ $^b$ Operations Research Center, Massachusetts Institute of Technology
 }
	
	
	%\toctitle{Lecture Notes in Computer Science}
	%\tocauthor{Authors' Instructions}
	\maketitle

\begin{abstract}
	\vspace{-1em}
The classical Canonical Correlation Analysis (CCA) identifies the correlations between two sets of multivariate variables based on their
covariance, which has been widely applied in diverse fields such as computer vision, natural language processing, and speech analysis. Despite its popularity, CCA can encounter challenges in explaining correlations between two variable sets within high-dimensional data contexts. Thus, this paper studies Sparse Canonical Correlation Analysis (SCCA) that enhances the interpretability of CCA. We first show that SCCA generalizes three well-known sparse optimization problems, sparse PCA, sparse SVD, and sparse regression, which are all classified as NP-hard problems. This result motivates us to develop strong formulations and efficient algorithms. Our main contributions include (i) the introduction of a combinatorial formulation that captures the essence of SCCA and allows the development of approximation algorithms; (ii) the derivation of an equivalent mixed-integer semidefinite programming model that facilitates a specialized branch-and-cut algorithm with analytical cuts; and (iii) the establishment of the complexity results for {two} low-rank special cases of SCCA. The effectiveness of our proposed formulations and algorithms is validated through numerical experiments.
%This paper studies the  Sparse Canonical Correlation Analysis (SCCA), which improves the interpretability of the classic CCA,  a powerful statistical technique used for multivariate data analysis. The fact that the SCCA is NP-hard motivates us to develop efficient formulations and algorithms. The main contributions of the paper are: (i) to develop an equivalent combinatorial formulation of SCCA; (ii) to derive an equivalent mixed-integer semidefinite program, which allows us design a customized branch-and-cut algorithm with the closed-form cuts; and (iii) to demonstrate that low-rank SCCA yields reduced complexity. Finally, we numerically test the proposed formulations and algorithms.
\end{abstract}

\vspace{-10pt}

\section{Introduction}
The Canonical Correlation Analysis (CCA), proposed by H. Hotelling  \cite{hotelling1935most}, aims to identify the correlations between two sets of multivariate variables based on their covariance. Since then, CCA has become a powerful statistical technique used for multivariate data analysis, with its applications across diverse fields such as computer vision \cite{huang2010super}, natural language processing \cite{vinokourov2002inferring}, and speech analysis \cite{hermansky1994rasta}. %However,  it is recognized that the CCA may have difficulty explaining the correlation between two sets of variables when being applied to high-dimensional datasets (e.g., genomic data \cite{parkhomenko2007genome}).Compared to CCA,  the Sparse CCA (SCCA) that seeks the sparse linear combinations of the two sets of variables allows for significantly better interpretability (see, e.g., \cite{witten2009penalized,witten2009extensions,yang2019survey}).
Despite its popularity, CCA can encounter challenges in explaining correlations between two variable sets within high-dimensional data contexts, such as genomic datasets \cite{parkhomenko2007genome}. In contrast, Sparse Canonical Correlation Analysis (SCCA), which seeks sparse linear combinations of these variable sets, offers substantially enhanced interpretability \cite{witten2009penalized,witten2009extensions,yang2019survey}.

Formally,
this paper studies the SCCA problem:
\begin{equation} \label{scca}
v^*:=	\max_{\bm x \in \Re^n, \bm y \in \Re^m} \left\{ \bm x^{\top} \bm A \bm y: \bm x^{\top}\bm B \bm x \le 1, \bm y^{\top}\bm C \bm y \le 1, \|\bm x\|_0 \le s_1, \|\bm y\|_0 \le s_2\right\}, \tag{SCCA}
\end{equation}
where $s_1\le n$, $s_2\le m$ are positive integers and $\begin{pmatrix}
    \bm B & \bm A\\
    \bm A^{\top} & \bm C
\end{pmatrix} $ denotes a covariance matrix of $(n+m)$ random variables. Specifically, $\bm B$ and $\bm C$ are the covariance matrices of the $n$ and $m$ random variables, respectively, and $\bm A\in \Re^{n\times m}$ is  the cross-covariance matrix between $n$ and $m$ random variables. Hence, $\begin{pmatrix}
    \bm B & \bm A\\
    \bm A^{\top} & \bm C
\end{pmatrix} $, $\bm B$, $\bm C$ are positive semidefinite matrices of size $(n+m), n$, and $m$, respectively.
Here, matrices $\bm B$, $\bm C$ can be singular, i.e., some random variables may be dependent on others. In fact, the covariance matrices $\bm B, \bm C$ are often low-rank, especially within the high-dimension low-sample size data context (see, e.g., the gene expression data in \cite{witten2009penalized}). 

The \ref{scca} problem 
%has a quite general framework
generalizes three widely-studied sparsity-constrained optimization problems as special cases, which are sparse PCA \cite{bertsimas2022solving,dey2022using,li2020exact}, sparse SVD \cite{li2021beyond,witten2009penalized}, and sparse regression \cite{hocking1967selection,bertsimasbest2016}. To be specific, when $n=m$, $s_1=s_2$, $\bm B, \bm C$ are identity matrices, and $\bm A$ is a positive semidefinite matrix, \ref{scca} reduces to the classic sparse PCA problem; when $\bm B, \bm C$ are identity matrices,  \ref{scca} becomes the sparse SVD problem; and when $\bm{A}$ is rank-one, \Cref{sec:rank} shows that \ref{scca} is equivalent to two sparse linear regression subproblems.

\vspace{-10pt}

\subsection{Main contributions}
\ref{scca} is generally NP-hard, given that its special cases, sparse PCA, sparse SVD, and sparse regression are all classified as NP-hard problems. We are motivated to develop efficient formulations and algorithms for \ref{scca} through a mixed-integer optimization lens. The main contributions, along with the structure of the remainder of this paper, are the following:
\begin{enumerate}[(i)]
    \item In \Cref{sec:com}, we present an exact semidefinite programming (SDP) reformulation and derive a closed-form optimal value of classic CCA problem. We also develop an equivalent combinatorial formulation of \ref{scca}; 
\item \Cref{sec:misdp} derives an equivalent mixed-integer SDP (MISDP) reformulation for \ref{scca}. When applying the Benders decomposition approach, instead of solving the large-scale SDPs, we design a customized branch-and-cut algorithm with closed-form cuts, which can successfully solve \ref{scca}  to optimality; 
\item When the covariance matrix $\begin{pmatrix}
    \bm B & \bm A\\
    \bm A^{\top} & \bm C
\end{pmatrix} $  is low-rank, \Cref{sec:rank} studies the complexity of two special cases of \ref{scca}; and
%\item By relaxing the binary variables in MISDP, we obtain an SDP relaxation used for providing an upper bound of \ref{scca};  and 
\item \Cref{sec:numerical} numerically test the proposed formulations and algorithms.
\end{enumerate}

\vspace{-15pt}

\subsection{Relevant literature}
\textit{SCCA.} To the best of our knowledge, the work \cite{parkhomenko2007genome} was the first paper that introduced the concept of SCCA to select only small subsets of variables to better explain the relationship between many genetic loci and gene expression phenotypes. A handful subset of features enhances interpretability, a desirable property, especially in complex data analysis, which has been successfully demonstrated in Sparse PCA \cite{jeffers1967two}. To obtain sparse canonical loadings $(\bm x, \bm y)$, \cite{ waaijenborg2008quantifying}  first applied elastic net penalty to the classical CCA via an iterative regression procedure. In a seminal work on SCCA \cite{witten2009penalized}, the authors proposed a rigorous formulation by enforcing the $\ell_1$ constraints on variables $(\bm x, \bm y)$ and developed a penalized matrix decomposition method to solve the penalized CCA problem.
Then, extensive research has focused on various penalty norm functions to obtain sparse canonical loadings (see, e.g.,  \cite{hardoon2011sparse,le2009sparse,  waaijenborg2008quantifying,witten2009extensions,chu2013sparse}). In particular, \cite{chu2013sparse} penalized multiple canonical loadings by $\ell_1$ norm and computed the sparse solution by the linearized Bregman method. It should be noted that under the assumption that the leading canonical loadings are sparse,  \cite{chen2013sparse,gao2015minimax,gao2017sparse} established theoretical guarantees of iterative approaches for estimating sparse solutions.
Another research direction in SCCA introduced penalty functions based on group structural information of input data and developed group SCCA methods \cite{lin2014correspondence,lin2013group}. For a comprehensive overview of CCA and SCCA methods, we refer readers to the survey by \cite{yang2019survey}  and the references therein. These approaches, however, do not strictly enforce the exact sparsity requirement but only approximate the sparsity requirement (i.e., the $\ell_0$ norm) by a convex function. Another relevant work \cite{watanabe2023branch} introduced binary variables to recast SCCA as a mixed-integer nonconvex program under the assumption of positive definite matrices $\bm B, \bm C$, based on which they designed a branch-and-bound algorithm. Different from the literature, our work does not require positive definiteness assumption of matrices $\bm B, \bm C$ and we are able to obtain mixed-integer conic and semidefinite programming reformulations, allowing for better exact and approximation algorithms.
%\vspace{1em}

\noindent\textit{Connections to and differences with sparse PCA and sparse SVD.} 
Analogous to \ref{scca}, both sparse PCA \cite{dey2022using,jeffers1967two} and sparse SVD \cite{li2021beyond} select small subsets of variables to improve the interpretability of dimensionality reduction methods: PCA and SVD.
Considerable investigation has been conducted on solving sparse PCA and sparse SVD from three angles: convex relaxations \cite{d2004direct,dey2022using,dey2023solving}, approximation algorithms \cite{chan2016approximability,chowdhury2020approximation,li2021beyond}, and exact algorithms \cite{bertsimas2022solving, li2020exact,li2021beyond}.
As mentioned before, in sparse PCA and sparse SVD, the covariance matrices $\bm B, \bm C$ are identity. Such a setting dramatically simplifies the subset selection problems of sparse PCA and sparse SVD  compared to that of SCCA, as in these problems, it suffices to focus on the selection of a submatrix of the matrix $\bm A$. Specifically, it is shown in \cite{d2008optimal, li2020exact, papailiopoulos2013sparse} that sparse PCA reduces to selecting a principal submatrix of $\bm A$ to maximize the largest eigenvalue(s) and sparse SVD reduces to selecting a possibly non-symmetric submatrix of $\bm A$ to maximize the largest singular value(s) \cite{li2021beyond}. Quite differently, the combinatorial reformulation \eqref{scca:com} of \ref{scca} aims to simultaneously select a sized-$(s_1\times s_1)$ principal submatrix of $\bm B$,  a sized-$(s_2\times s_2)$ principal submatrix of $\bm C$, and a sized-$(s_1\times s_2)$ submatrix of $\bm A$. These fundamental differences in the underlying formulations of sparse PCA and sparse SVD preclude the direct application of their existing algorithms to the \ref{scca}. %This fundamental difference prevents the existing formulations and algorithms developed for either sparse PCA or sparse SVD  from directly solving \ref{scca}.

\noindent\textbf{Notations:} 
The following notation is used throughout the paper. We use bold lower-case letters (e.g., $\bm{x}$) and bold upper-case letters (e.g., $\bm{X}$) to denote vectors and matrices, respectively, and we use corresponding non-bold letters (e.g., $x_i$) to denote their components. We let $\S^n,\S_+^n, \S_{++}^n$ denote the set of all the $n\times n$ symmetric real matrices, the set of all the $n\times n$ symmetric positive semidefinite matrices, and the set of all the $n\times n$ symmetric positive definite matrices, respectively. 
We let $\bm I$ denote the identity matrix and let $\bm 0$ denote the vector or matrix with all-zero entries.
We let $\Re^n_+$ denote the set of all $n$-dimensional nonnegative vectors. We let $[n]:=\{1,2,\cdots, n\}$, $[s,n]:=\{s,s+1,\cdots, n\}$. Given a matrix $\bm{A}\in \Re^{n\times m}$ and two subsets $S\subseteq [n]$, $T\subseteq [m]$, we let $\bm A^{\dag}$ denote the pseudo inverse of matrix $\bm A$, let $\bm{A}_{S,T}$ denote a submatrix of $\bm{A}$ with rows and columns indexed by sets $S, T$, respectively, and let $(\bm{A}_{S,T})^{\dag}$ denote the pseudo inverse of submatrix $\bm{A}_{S,T}$. For a set $S$ and an integer $k$, we define the set $S+k:=\{i+k| i\in S\}$.
Given a vector $\bm{a}\in \Re^{n}$ and a subset $S\subseteq [n]$, we  let $\bm{a}_{S}$ denote a subvector of $\bm{a}$ in the subset $S$. We define $[\lambda]_+:=\max\{\lambda, 0\}$.
We let $\sigma_{\max}(\cdot)$ denote the largest singular value function and let $\lambda_{\max}(\cdot)$ denote the largest eigenvalue value function.

\vspace{-10pt}

\section{A combinatorial reformulation of SCCA} \label{sec:com}
%This section proposes an equivalent combinatorial optimization reformulation of \ref{scca}, which allows us to design two effective approximation algorithms to solve it.
% as selecting two sized-$s_1$, $s_2$ subsets from $[n], [m]$ to maxim
This section introduces an equivalent combinatorial optimization reformulation of \ref{scca}. This reformulation serves as the foundation for developing two effective approximation algorithms to solve \ref{scca} in \Cref{sec:numerical}.

\vspace{-10pt}

\subsection{An exact semidefinite programming  representation of CCA}
To begin with, let us focus on the classic CCA problem, which refers to \ref{scca} without zero-norm constraints, as defined below:
\begin{equation} \label{cca}
\max_{\bm x \in \Re^n, \bm y \in \Re^m} \left\{ \bm x^{\top} \bm A \bm y: \bm x^{\top}\bm B \bm x \le 1, \bm y^{\top}\bm C \bm y \le 1\right\}. \tag{CCA}
\end{equation}
This formulation of \ref{cca} can be regarded as a quadratically constrained quadratic program concerning the variables $\begin{pmatrix}
    \bm x\\
    \bm y
\end{pmatrix} \in \Re^{n\times m}$. 
% That is, the \ref{cca} can be equivalently converted into
% \begin{equation} \label{cca1}
% \max_{\bm x \in \Re^n, \bm y \in \Re^m} \left\{ \begin{pmatrix}
%     \bm x\\
%     \bm y
% \end{pmatrix}^{\top} \tilde{\bm A} \begin{pmatrix}
%     \bm x\\
%     \bm y
% \end{pmatrix}: \begin{pmatrix}
%     \bm x\\
%     \bm y
% \end{pmatrix}^{\top} \tilde{\bm B} \begin{pmatrix}
%     \bm x\\
%     \bm y
% \end{pmatrix} \le 1, \begin{pmatrix}
%     \bm x\\
%     \bm y
% \end{pmatrix}^{\top} \tilde{\bm C} \begin{pmatrix}
%     \bm x\\
%     \bm y
% \end{pmatrix} \le 1\right\}, 
% \end{equation}
We next define three-block matrices of size $(n+m)$ below that aid in the presentation of our results.
\begin{align*}
\tilde{\bm A}=\begin{pmatrix}
\bm 0 & \bm A/2\\
\bm A^{\top}/2 & \bm 0
\end{pmatrix}, \ \ \tilde{\bm B}=\begin{pmatrix}
\bm B & \bm 0\\
\bm 0& \bm 0
\end{pmatrix}, \ \ \tilde{\bm C}=\begin{pmatrix}
\bm 0 & \bm 0\\
\bm 0& \bm C
\end{pmatrix}.
\end{align*}

By introducing a size-$(n+m)$ matrix variable $\bm X = \begin{pmatrix}
    \bm x\\
    \bm y
\end{pmatrix} \begin{pmatrix}
    \bm x\\
    \bm y
\end{pmatrix}^{\top}$ and removing the rank-one constraint on $\bm X$, we can obtain an SDP relaxation of \eqref{cca}, as described below
\begin{equation}\label{sdp}
	\max_{\bm X \in \S_+^{m+n}}  \left\{ \tr\left(\tilde{\bm A} \bm X\right): \tr\left( \tilde{\bm B} \bm X\right) \le 1, \tr\left(\tilde{\bm C}\bm X\right)\le 1\right\}. \tag{SDP Relaxation}
\end{equation}

Next, let us present a key lemma regarding properties of block matrices being positive semidefinite, fundamental for reformulating the \ref{scca}.
\begin{lemma}[\cite{gallier2020schur}]\label{lem:psd}
    For any symmetric block matrix $\begin{pmatrix}
    \bm B & \bm A\\
    \bm A^{\top} & \bm C
\end{pmatrix} \in \S^{n+m}$, the followings are equivalent:
\begin{enumerate}[(i)]
\item The block matrix is positive semidefinite; 
\item $\bm B \in \S_+^n$, $(\bm I-\bm B \bm B^{\dag})\bm A =\bm 0$, $\bm C -\bm A^{\top}\bm B^{\dag} \bm A \in \S_+^m$; and 
\item $\bm C\in \S_+^m$, $(\bm I-\bm C \bm C^{\dag})\bm A ^{\top}=\bm 0$, $\bm B -\bm A\bm C^{\dag} \bm A^{\top} \in \S_+^n$.
\end{enumerate}
\end{lemma}
%According to \Cref{lem:psd}, we show that the classic CCA problem can be equivalently converted into a semidefinite program. 

Inspired by \Cref{lem:psd}, we hereby establish the equivalence between \ref{cca} and its \ref{sdp}. Remarkably, both of these problems achieve the same optimal value, namely $\sigma_{\max}(\sqrt{\bm B^{\dag}} \bm A \sqrt{\bm C^{\dag}})$. %we demonstrate below that \ref{cca} and its \ref{sdp} are equivalent, both attaining the optimal value $\sigma_{\max}(\sqrt{\bm B^{\dag}} \bm A \sqrt{\bm C^{\dag}})$.

\begin{proposition}\label{prop:cca}
For the \ref{cca}, we have the following results.
\begin{enumerate}[(i)]
    \item Both \ref{cca} and its \ref{sdp} have an optimal value $\sigma_{\max}(\sqrt{\bm B^{\dag}} \bm A \sqrt{\bm C^{\dag}})$;
    \item A pair of optimal solutions $(\bm x^*, \bm y^*)$ to \ref{cca} satisfies 
    $$\bm x^* = \sqrt{\bm B^{\dag}} \bm q, \ \ \bm y^* = \sqrt{\bm C^{\dag}} \bm p,$$ 
    where $\bm q \in \Re^n, \bm p \in \Re^m$ denote a pair of leading singular vectors of matrix $\sqrt{\bm B^{\dag}} \bm A \sqrt{\bm C^{\dag}}$; and
    \item  An optimal solution $\bm X^*$ to the \ref{sdp} is
    $$\bm X^* =  \begin{pmatrix}
    \bm x^*\\
    \bm y^*
\end{pmatrix} \begin{pmatrix}
    \bm x^*\\
    \bm y^*
\end{pmatrix}^{\top}.$$
\end{enumerate}
\end{proposition}
 \noindent\textit{Proof.}
    See Appendix A.1. \qed

The proof of \Cref{prop:cca} motivates the following observation on the optimal values of \ref{cca} and \ref{scca}.
\begin{observation}\label{remark1}
  The optimal value of  \ref{cca} is upper bounded by  $1$, so is the optimal value of  \ref{scca}.
\end{observation}
 \noindent\textit{Proof.}
Since matrix $\begin{pmatrix}
    \bm B & {\bm A}\\
{\bm A^{\top}}&  \bm C
\end{pmatrix}$ denotes a covariance matrix of a subset of variables and thus is always positive semidefinite. According to \Cref{lem:psd}, we have that
$$\bm B \succeq \bm A \bm C^{\dag} \bm A^{\top} \Longrightarrow \bm I \succeq \sqrt{\bm B^{\dag}}  \bm A  \bm C^{\dag} \bm A^{\top} \sqrt{\bm B^{\dag}} ,$$
which means that $\sigma_{\max}\left(\sqrt{\bm B^{\dag}} \bm A \sqrt{\bm C^{\dag}}\right)\le 1$ must hold.  \qed

It is noteworthy that the results presented in Proposition \ref{prop:cca} are established through a distinct methodology. This methodology leverages the positive semidefinite condition of block matrices, as shown in Lemma \ref{lem:psd}, and incorporates duality theory. This approach differs from most prior research \cite{lu2014large,parkhomenko2009sparse,yang2019survey}, which proved Part (i) of Proposition \ref{prop:cca} by relying on the singular value decomposition and assuming that matrices $\bm B$ and $\bm C$ are positive definite (i.e., full rank). To the best of our knowledge, \cite{chu2013sparse}  showed parts (i) and (ii) of Proposition \ref{prop:cca} for a special low-rank \ref{cca} problem, where the authors assumed that the covariance matrices are defined as $\bm A = \bm U \bm V^{\top}$,  $\bm B = \bm U \bm U^{\top}$,  and $\bm C = \bm V \bm V^{\top}$. 
Remarkably, Proposition \ref{prop:cca} extends this result to a more general scenario where $\bm B$ and $\bm C$ are not constrained to be strictly positive definite and $\bm A$ is not constrained to directly depend on $\bm B, \bm C$, allowing for rank deficiencies and flexible data structure.

%the results in \Cref{prop:cca} are established by leveraging the positive semidefinite condition of block matrices  in \Cref{lem:psd} and duality theory. This differs from the existing work in that they proved Part (i) of \Cref{prop:cca}   based on the singular value decomposition, under the assumption that matrices $\bm B, \bm C$ are positive definite, i.e., full rank (see, e.g., \cite{lu2014large,parkhomenko2009sparse,yang2019survey}). \Cref{prop:cca} shows that the result still holds in the general rank-deficient case, i.e.,   $\bm B, \bm C$ are not required to be strictly positive definite.

\vspace{-10pt}

\subsection{An equivalent formulation of \ref{scca}}
%This subsection reduces \ref{scca} to the combinatorial optimization based on the results in \Cref{prop:cca}.
In this subsection, we transform \ref{scca} into a combinatorial optimization problem, according to the insights provided by Proposition \ref{prop:cca}.
\begin{theorem}\label{them:com}
The \ref{scca} is equivalent to the following combinatorial optimization:
\begin{align}\label{scca:com}
v^*:=\max_{\begin{subarray}{c}
			S_1\subseteq [m], |S_1|\le s_1, \\
			S_2 \subseteq [n], |S_2|\le s_2
			\end{subarray}}\left\{\sigma_{\max}\left(\sqrt{(\bm B_{S_1, S_1})^{\dag}} \bm A_{S_1, S_2} \sqrt{(\bm C_{S_2, S_2})^{\dag}}\right)\right\}.
\end{align}
\end{theorem}
 \noindent\textit{Proof.}
By introducing the  subsets $(S_1, S_2)$ to denote the supports of variables $(\bm x, \bm y)$ in \ref{scca},  then we can remove the zero-norm constraints on $(\bm x, \bm y)$ and reformulate \ref{scca} as
 \begin{align} \label{scca_set}
     &v^*:= \max_{\begin{subarray}{c}
			S_1\subseteq [m], |S_1|\le s_1, \\
			S_2 \subseteq [n], |S_2|\le s_2
			\end{subarray}} \max_{
   \begin{subarray}{c}
   \bm x \in \Re^{|S_1|},\\ \bm y \in \Re^{|S_2|}
   \end{subarray}} \left\{ \bm x^{\top} \bm A_{S_1, S_2} \bm y:  \bm x^{\top}\bm B_{S_1, S_1} \bm x \le 1, \bm y^{\top}\bm C_{S_2, S_2} \bm y \le 1\right\}. 
 \end{align}

Following from the Part (i) in \Cref{prop:cca}, we can show that
 for any subsets $S_1\subseteq [n], S_2\subseteq [m]$, the following identity holds.
      \begin{align*}
    & \max_{   \begin{subarray}{c}
   \bm x \in \Re^{|S_1|}, \bm y \in \Re^{|S_2|}
   \end{subarray}} \left\{ \bm x^{\top} \bm A_{S_1, S_2} \bm y:  \bm x^{\top}\bm B_{S_1, S_1} \bm x \le 1, \bm y^{\top}\bm C_{S_2, S_2} \bm y \le 1\right\} \\
   &= \sigma_{\max}\left(\sqrt{(\bm B_{S_1, S_1})^{\dag}} \bm A_{S_1, S_2} \sqrt{(\bm C_{S_2, S_2})^{\dag}}\right).
 \end{align*}
Plugging the result above into the inner maximization problem in \eqref{scca_set}, we complete the proof.
\qed

The combinatorial formulation \eqref{scca:com} 
presents significant computational difficulties when attempting to solve \ref{scca}. The primary obstacles are two-fold: first, simultaneously selecting submatrices from the matrices $\bm A, \bm B, \bm C$ requires a sophisticated optimization across multiple dimensions. Second, the selection criterion is particularly complex, as it involves optimizing the largest singular value of the product of the selected submatrix of $\bm A$ and the square root of pseudo-inverse submatrices of $\bm B$ and $\bm C$. These complexities necessitate effective optimization solution procedures to address the high-dimensional and non-convex nature of the problem. %poses challenges of solving \ref{scca} in two aspects: (i)  simultaneous submatrix selections of $\bm A, \bm B, \bm C$ and (ii) a complex selection criterion involving the largest singular value and the square root of pseudo-inverse matrices.

As a side product of \Cref{remark1}, the optimal value of  \ref{scca} is trivially upper bounded by 1.
\begin{observation}
  The optimal value of  \ref{scca}  satisfies  $v^*\le 1$.
\end{observation}

\vspace{-10pt}

\section{Reformulating \ref{scca} as a mixed-integer semidefinite program (MISDP)}\label{sec:misdp}
%This section develops an equivalent MISDP formulation of \ref{scca} that motivates the design of a branch-and-cut algorithm.
This section formulates an equivalent Mixed-Integer Semidefinite Programming (MISDP) formulation for the \ref{scca} problem. This reformulation serves as the foundation for developing a branch-and-cut algorithm to solve the problem effectively.
\subsection{Valid inequalities for \ref{scca}}
We prove that there exists a bounded optimal solution $(\bm x^*, \bm y^*)$ of the \ref{scca}. To be specific, we show that there exists an optimal solution $(\bm x^*, \bm y^*)$ of the \ref{scca} satisfying the constraints $\|\bm x^*\|^2_2\le M_1$ and $\|\bm y^*\|^2_2\le M_2$, where $M_1$ and $M_2$ are finite-valued parameters. 

% Therefore, we can equivalently recast  \ref{scca} into
% \begin{equation} \label{scca:bound}
% \begin{aligned}
% v^*:=	\max_{\bm x \in \Re^n, \bm y \in \Re^m} \big\{ \bm x^{\top} \bm A \bm y: & \bm x^{\top}\bm B \bm x \le 1, \|\bm x\|_0 \le s_1, \|\bm x\|^2_2\le M_1,\\ 
% &\bm y^{\top}\bm C \bm y \le 1, \|\bm y\|_0 \le s_2, \|\bm y\|^2_2\le M_2\big\}.
% \end{aligned}
% \end{equation}

\begin{proposition}\label{prop:bigm}
    The \ref{scca} admits an optimal solution $(\bm x^*, \bm y^*)$ satisfying $\|\bm x^*\|^2_2\le M_1$ and $\|\bm y^*\|^2_2\le M_2$, where $M_1:={1/\lambda_r(\bm{B})+1/(\lambda_r(\bm{B})s_{\min}(\bm{B}))}$ and $M_2:={1/\lambda_{\hat{r}}(\bm{C})+1/(\lambda_{\hat{r}}(\bm{C})s_{\min}(\bm{C}))}$ with $\lambda_r(\bm{B}),\lambda_{\hat{r}}(\bm C)$ being the smallest nonzero eigenvalues of matrices $\bm B,\bm C$ and $s_{\min}(\bm{R})$ being the smallest nonzero singular value of all the submatrices of the zero eigenvectors of matrix $\bm{R}$.
\end{proposition}
 \noindent\textit{Proof.}
    See Appendix A.2. \qed

The proof of Proposition~\ref{prop:bigm} is straightforward in the case when $\bm B$ and $\bm C$ are of full rank as in this case the feasible region is a bounded set. In order to prove the result in the case when $\bm B$ is not full-rank, one has to show that it is possible to construct sparse solutions that are not ``too far" away.

In fact, the bounds $M_1,M_2$ in \Cref{prop:bigm} also hold for any given feasible subsets $(S_1,S_2)$ of SCCA \eqref{scca:com}.
\begin{corollary}\label{cor_prop:bigm}
For any given feasible subsets $(S_1,S_2)$ of SCCA \ref{scca:com}, there exists a SCCA feasible solution $(\bm x, \bm y)$ such that the supports of $\bm x, \bm y$ are $S_1,S_2$, respectively and we have that $\|\bm x\|^2_2\le M_1$ and $\|\bm y\|^2_2\le M_2$, where $M_1,M_2$ are defined in \Cref{prop:bigm}.
\end{corollary}

\vspace{-10pt}

\subsection{An equivalent MISDP formulation}
While the combinatorial formulation \eqref{scca:com} is elegant in its structure, it poses significant challenges when attempting to solve it to optimality using branch-and-bound based methods. To fill this gap, in this subsection, we derive an equivalent MISDP formulation for \ref{scca}, amenable for developing exact methods.

It is convenient to define the following notation. Let $M_{ii}$ be defined as follows:
\begin{align*}
    M_{ii} = \begin{cases}
        M_1, &\forall i \in [n],\\
M_2, & \forall i \in [n+1,n+m].
    \end{cases}
\end{align*}
\begin{theorem}\label{them:misdp}
The  \ref{scca} is equivalent to the following MISDP: 
\begin{equation}\label{eq:misdp}
	v^*:=	\max_{\begin{subarray}{c}  
	\bm X \in \S_+^{n+m}, \\\bm z \in \mathcal{Z}\end{subarray}} \{\tr(\tilde{\bm A} \bm X): \tr( \tilde{\bm B} \bm X) \le 1, 
 \tr(\tilde{\bm C} \bm X)\le 1, X_{ii} \le M_{ii} z_i, \forall i\in [n+m] \}.
\end{equation}
where the feasible set of variables $\bm z$ is defined as $\Z:=\{\bm z \in \{0,1\}^{n+m}:  \sum_{i\in [n]} z_i \le s_1,  \sum_{i\in [n+1, n+m]} z_i \le s_2\}$. 
%and throughout this paper, 
%the coefficient $M_{ii}$ is defined as
% \begin{align*}
%     M_{ii} = \begin{cases}
%         M_1, &\forall i \in [n];\\
% M_2, & \forall i \in [n+1,n+m].
%     \end{cases}
% \end{align*}
% \begin{align*}
%     M_{ij} = \begin{cases}
%         M_1, &\forall i \in [n], j\in [n], i=j;\\
%     M_1/2, & \forall i \in [n], j\in [n], i\neq j;\\
% M_2, & \forall i \in [n+1,n+m], j\in [n+1, n+m], i= j;\\
%     M_2/2, & \forall i \in [n+1,n+m], j\in [n+1, n+m], i\neq j;\\
%        M_1/2+  M_2/2,  & i \in [n], j\in [n+1, n+m] \text{ or } \forall i \in [n+1,n+m], j\in [n].
%     \end{cases}
% \end{align*}
\end{theorem}
\noindent\textit{Proof.}
For the SCCA \eqref{scca_set}, according to \Cref{prop:cca}, the inner maximization problem admits an exact semidefinite programming formulation. Using the  variables $\bm z\in \Z$ to describe the set constraints in  SCCA \eqref{scca_set}, we can reformulate it as
\begin{equation} \label{scca:binary}
\begin{aligned}
v^*:=\max_{\bm z\in \Z}	\max_{\bm X \in  \S_+^{n+m}} \big\{\tr(\tilde{\bm A} \bm X): &\tr( \tilde{\bm B} \bm X) \le 1, 
 \tr(\tilde{\bm C} \bm X)\le 1,\\
 &X_{ii} (1-z_i) = 0, \forall i\in [m+n]\big\}.
\end{aligned}
\end{equation}

\Cref{prop:bigm} shows that there is an optimal solution $(\bm x^*, \bm y^*)$ to \ref{scca} that satisfies $\|\bm x^*\|_2^2 \le M_1$ and  $\|\bm y^*\|_2^2 \le M_2$. Based on this, we can construct an optimal solution $(\bm z^*, \bm X^*)$ 
for SCCA \eqref{scca:binary} by letting
\begin{align*}
\bm X^* = \begin{pmatrix}
    \bm x^*\\
    \bm y^*
\end{pmatrix} \begin{pmatrix}
    \bm x^*\\
    \bm y^*
\end{pmatrix}^{\top},
z_i = \begin{cases}
    1 & \text{if } x_i^* \neq 0\\
    0 & \text{if } x_i^* = 0
\end{cases}, \forall i\in [n], z_{i+n} = \begin{cases}
    1 & \text{if } y_i^* \neq 0\\
    0 & \text{if } y_i^* = 0
\end{cases}, \forall i\in [m],
\end{align*}
where the optimal solution $\bm X^*$ satisfies the following inequalities
\begin{align*}
 &X_{ii}^* = (x_i^*)^2 \le M_1 z_i, \forall i\in [n],\ \  X_{(i+n) (i+n)}^* = (y_i^*)^2 \le M_2 z_{i+n}, \forall i\in [m].
%  &|X_{ij}^*| = x_i^* x_j^* \le \frac{(x_i^*)^2 + (x_j^*)^2}{2} \le \frac{M_1}{2}z_i, \forall i\in [n],j \in [n], i\neq j,\\
%   &|X_{(i+n)(j+n)}^*| = y_i^* y_j^* \le \frac{(y_i^*)^2 + (y_j^*)^2}{2} \le \frac{M_2}{2}z_{i+n}, \forall i\in [m],j \in [m], i\neq j,\\
% &|X_{(i)(j+n)}^*| = x_i^* y_j^* \le \frac{(x_i^*)^2 + (y_j^*)^2}{2} \le (\frac{M_1}{2} + \frac{M_2}{2}) z_i, \forall i\in [n],j \in [m], i\neq j,
\end{align*}
This allows us to recast the SCCA \eqref{scca:binary} into an MISDP formulation \eqref{eq:misdp}.
\qed

Note that the proposed MISDP formulation \eqref{eq:misdp} is of size $(n+m)\times (n+m)$ since our matrix variable $\bm{X}$ replaces $\begin{pmatrix}
    \bm x\\
    \bm y
\end{pmatrix} \begin{pmatrix}
    \bm x\\
    \bm y
\end{pmatrix}^{\top}$ in \ref{scca}.
Relaxing the binary variables in SCCA \eqref{eq:misdp}  to be continuous, we obtain an upper bound of SCCA \eqref{eq:misdp}, i.e., $	v^* \le \hat v$
\begin{equation}\label{eq:sdp}
\hat v:=\max_{\bm X \in \S_+^{n+m}, \bm z \in \hat{Z}} \{\tr(\tilde{\bm A} \bm X): \tr( \tilde{\bm B} \bm X) \le 1, 
 \tr(\tilde{\bm C} \bm X)\le 1, X_{ii} \le M_{ii} z_i, \forall i\in [n+m]\}.
\end{equation}
where  $\hat{\Z}:=\{\bm z \in [0,1]^{n+m}:  \sum_{i\in [n]} z_i \le s_1,  \sum_{i\in [n+1, n+m]} z_i \le s_2\}$. This SDP relaxation \eqref{eq:sdp} can be directly solved by commercial solvers such as MOSEK or SDPT3.

\vspace{-10pt}

\subsection{Developing a branch-and-cut algorithm with closed-form cuts}
By dualizing the inner maximization problem over $\bm X$ in the  MISDP \eqref{eq:misdp}, in this subsection, we derive an equivalent mixed-integer linear program for \ref{scca}, which motivates us to develop a branch-and-cut algorithm. 

By separating the binary variables $\bm z$, we rewrite the MISDP \eqref{eq:misdp} as
\begin{align}\label{scca_z}
v^*:=	\max_{\bm z \in \Z, v} \{v: v\le f(\bm z)\},
\end{align}
where  the  function $f(\bm z)$ is defined as
\begin{equation}\label{eq:obj}
	f(\bm z):=	\max_{\bm X \in \S_+^{n+m}} \left\{\tr(\tilde{\bm A} \bm X): \tr( \tilde{\bm B} \bm X) \le 1, 
 \tr(\tilde{\bm C} \bm X)\le 1, X_{ii} \le M_{ii} z_i, \forall i\in [n+m]\right\}.
\end{equation}
By introducing the Lagrangian multipliers $(\theta_1, \theta_2, \bm \lambda)$, the Lagrangian dual of the maximization problem \eqref{eq:obj} can be written as
\begin{equation}\label{eq:dual}
	\begin{aligned}
	f(\bm z)&= \min_{\begin{subarray}{c}
	     \theta_1 \ge 0, \theta_2 \ge 0, \\ \bm \lambda \in \Re_+^{n+m}
	\end{subarray}}	\max_{\bm X \in \S_+^{n+m}} \tr(\tilde{\bm A} \bm X) -\theta_1 \tr( \tilde{\bm B} \bm X) -\theta_2 \tr( \tilde{\bm C} \bm X) + \theta_1 + \theta_2,\\
 & -\sum_{i\in [n+m]} \lambda_i X_{ii} + \sum_{i\in [n+m]} \lambda_i M_{ii}  z_i\\
 & = \min_{\begin{subarray}{c}
	     \theta_1 \ge 0, \theta_2 \ge 0, \\ \bm \lambda \in \Re_+^{n+m}
	\end{subarray}}	\bigg\{ \theta_1 + \theta_2 + \sum_{i\in [n+m]} \lambda_i M_{ii}  z_i: \begin{pmatrix}
	   & \theta_1\bm B  & -\bm A/2\\
     & -\bm A^{\top}/2 & \theta_2\bm C
	\end{pmatrix} \succeq - \Diag(\bm \lambda) \bigg\},
	\end{aligned} 
\end{equation}
where the strong duality holds due to the function $f(\bm{z})$ being concave, bounded, and thus continuous in the set $\hat{\Z}$ and Slater condition holds for any interior point $\bm{z}$ in the set $\hat{\Z}$.

Below, we derive the closed-form expression of the function $f(\bm z)$ with the given binary variable $\bm z\in \Z$. This allows us to reformulate SCCA \eqref{scca_z} as a mixed-integer linear program with exponentially many linear constraints and an efficient separation oracle.

\begin{proposition}\label{prop:milp}
The SCCA \eqref{scca_z} is  equivalent to 
\begin{equation} 
\begin{aligned}
v^*:=	&\max_{\bm z \in \Z, v} \bigg\{v: v\le  \sigma_{\max}\left(\sqrt{(\bm B_{S_1, S_1})^{\dag}} \bm A_{S_1, S_2} \sqrt{(\bm C_{S_2, S_2})^{\dag}}\right) +   \\
&\sum_{i\in S_1 \cup (S_2+n)} \lambda^* M_{ii}z_i: \forall S_1 \subseteq [n], |S_1|\le s_1, S_2 \subseteq [m], |S_2|\le s_2\bigg\},
\end{aligned}
\label{scca:milp}
\end{equation}
where for a pair of subsets $(S_1, S_2)$, the scalar $\lambda^*$ is defined as the largest positive eigenvalue of matrix $\bm D_2^{\top} \bm D_1^{-1} \bm D_2 -\bm D_3$ with $$ \bm D_1:= \begin{pmatrix}
	   & \theta^*_1\bm B_{S_1, S_1}  & -\bm A_{S_1, S_2}/2\\
     & -\bm A_{S_1, S_2}^{\top}/2 & \theta^*_2\bm C_{S_2, S_2}
	\end{pmatrix}, \ \ \bm D_2:=\begin{pmatrix}
	   & \theta^*_1\bm B_{S_1, [n]\setminus S_1}  & -\bm A_{S_1, [m]\setminus S_2}/2\\
     & -\bm A_{S_2, [n]\setminus S_1}^{\top}/2 & \theta^*_2\bm C_{S_2, [m]\setminus S_2}
	\end{pmatrix},$$
 and
 $$ \bm D_3:=   \begin{pmatrix}
	   & \theta^*_1\bm B_{[n]\setminus S_1, [n]\setminus S_1}    & -\bm A_{[n]\setminus S_1, [m]\setminus S_2}/2\\
     & -\bm A_{[n]\setminus S_1, [m]\setminus S_2}^{\top}/2 & \theta^*_2\bm C_{[m]\setminus S_2, [m]\setminus S_2}
	\end{pmatrix},$$
where $\theta^*_1 = \theta^*_2=\sigma_{\max}\left(\sqrt{(\bm B_{S_1, S_1})^{\dag}} \bm A_{S_1, S_2} \sqrt{(\bm C_{S_2, S_2})^{\dag}}\right)/2$.
\end{proposition}
 \noindent\textit{Proof.}
    See Appendix A.3. \qed

We note that SCCA \eqref{scca:milp} can be implemented via a delayed cut-generation procedure. That is, at each feasible branch-and-bound node with a binary solution $\hat{\bm z}$, let ${S}_1:=\{i: \hat{z}_i=1,\forall i\in [n]\}$ and ${S}_2:=\{i-n: \hat{z}_i=1,\forall i\in [n+1,n+m]\}$. Then we can compute the corresponding scalar $\lambda^*$ and generate the following valid inequality based on \eqref{scca:milp}:
\[v\le  \sigma_{\max}\left(\sqrt{(\bm B_{S_1, S_1})^{\dag}} \bm A_{S_1, S_2} \sqrt{(\bm C_{S_2, S_2})^{\dag}}\right) +  \sum_{i\in S_1 \cup (S_2+n)} \lambda^* M_{ii}z_i.\]

\vspace{-10pt}

\section{Low-rank SCCA}\label{sec:rank}
In practice, it is common that the sample covariance matrix $\begin{pmatrix}
\bm B & \bm A\\
\bm A^{\top} & \bm C
\end{pmatrix}$ exhibits low-rank characteristics. This phenomenon is especially prominent when dealing with high-dimensional, low-sample size data (e.g., gene expression data \cite{witten2009penalized}). In this section, we study two special cases of low-rank \ref{scca} and their computational complexities. %implications for solving the problem denoted as \ref{scca}.

% The  sample covariance matrix  $\begin{pmatrix}
%     \bm B & \bm A\\
%     \bm A^{\top} & \bm C
% \end{pmatrix} $  is often low-rank in practice, especially in the fields where high-dimension low-sample size data are prevalent (see, e.g., the gene expression data in \cite{witten2009penalized}). This section presents two low-rank cases and their effect on solving \ref{scca}. 

\subsection{Special Case I: \ref{scca} with low-rank covariance matrices}
In this section, we show that the computational complexity of \ref{scca} is contingent upon the ranks of the covariance matrices $\bm B$ and $\bm C$. To be more precise, when the sparsity level $s_1$ (or $s_2$) is equal to or greater than the rank of the covariance matrix $\bm B$ (or $\bm C$), the imposition of a zero-norm constraint over $\bm x$ (or $\bm y$) in \ref{scca} becomes redundant. Consequently, lower ranks in the covariance matrices correspond to better computational complexity in solving \ref{scca}.
%In this section, we show that the ranks of covariance matrices $\bm B, \bm C$  decide the computational complexity of \ref{scca}. More specifically, when the sparsity level $s_1$ (or $s_2$) is no less than the rank of covariance matrix $\bm B$ (or $\bm C$),  the zero-norm constraint over $\bm x$ (or $\bm y$) in \ref{scca} is redundant. As a result,  the lower rank implies a lower complexity of \ref{scca}. 
\begin{theorem}\label{them:rank}
Suppose $r:=\rank(\bm B)$ and $\hat r:=\rank(\bm C)$,  then the  \ref{scca}  takes a complexity of $\mathcal{O}(n^{r-1}m^{\hat r-1} + n^{r-1} + m^{\hat r-1})$. The following results hold:
\begin{enumerate}[(i)]
%     \item The \ref{cca} admits an optimal solution $(\bm x^*,\bm y^*)$ satisfying $\|\bm x^*\|_0\le r$ and  $\|\bm y^*\|_0\le \hat r$. That is, the \ref{cca}  is equivalent to 
%         \begin{equation} \label{scca_re}
% \max_{\bm x \in \Re^n, \bm y \in \Re^m} \left\{ \bm x^{\top} \bm A \bm y: \bm x^{\top}\bm B \bm x \le 1, \bm y^{\top}\bm C \bm y \le 1, \|\bm x\|_0 \le r, \|\bm y\|_0 \le \hat r\right\};
% \end{equation}
\item When $s_1\ge r$ and $s_2 \ge \hat r$,  the \ref{scca} problem is equivalent to \ref{cca}, i.e.,
        \begin{equation} \label{scca_re}
v^*:=\max_{\bm x \in \Re^n, \bm y \in \Re^m} \left\{ \bm x^{\top} \bm A \bm y: \bm x^{\top}\bm B \bm x \le 1, \bm y^{\top}\bm C \bm y \le 1\right\};
\end{equation}
\item When $s_1\ge r$ and $s_2 < \hat r$,  the \ref{scca} problem can be reduced to
\begin{equation} \label{scca_re2}
v^*:=\max_{\bm x \in \Re^n, \bm y \in \Re^m} \left\{ \bm x^{\top} \bm A \bm y: \bm x^{\top}\bm B \bm x \le 1, \bm y^{\top}\bm C \bm y \le 1, \|\bm y\|_0 \le s_2\right\};
\end{equation}
\item When $s_1< r$ and $s_2 \ge \hat r$,  the \ref{scca} problem can be reduced to
\begin{equation} \label{scca_re3}
v^*:=\max_{\bm x \in \Re^n, \bm y \in \Re^m} \left\{ \bm x^{\top} \bm A \bm y: \bm x^{\top}\bm B \bm x \le 1, \bm y^{\top}\bm C \bm y \le 1, \|\bm x\|_0 \le s_1\right\}.
\end{equation}
%     \item 
%     The \ref{scca} problem can be reduced to
%     \begin{equation}
%     \begin{aligned} \label{scca_re}
% v^*:=	\max_{\bm x \in \Re^n, \bm y \in \Re^m} \{ \bm x^{\top} \bm A \bm y: &\bm x^{\top}\bm B \bm x \le 1, \bm y^{\top}\bm C \bm y \le 1, \\
% &
% \|\bm x\|_0 \le \min\{r, s_1\}, \|\bm y\|_0 \le \min\{\hat r, s_2\}\}.
% \end{aligned}
%  \end{equation}
% Then, 
\end{enumerate}
\end{theorem}
 \noindent\textit{Proof.}
    See Appendix A.4. \qed

The proof of \Cref{them:rank} implies that \ref{cca} admits an optimal sparse solution $(\bm x^*, \bm y^*)$ satisfying  $\|\bm x^*\|_0\le r$ and  $\|\bm y^*\|_0\le \hat r$, provided that $\bm B, \bm C$ are of rank-$r$, $\hat r$, respectively. Thus, \Cref{them:rank} establishes a sufficient condition (i.e., $s_1\le r, s_2\le \hat r$) about when  \ref{cca} can be equivalent to \ref{scca}. Besides, \Cref{them:rank} implies the complexity of solving \ref{scca}, as summarized below. %It is worth noting that our proposed complexity matches that for low-rank sparse PCA in \cite{del2023sparse}.

\begin{corollary}
Suppose $r:=\rank(\bm B)$ and $\hat r:=\rank(\bm C)$. There exists an algorithm that can find an optimal solution to \ref{scca} in   $\mathcal{O}(n^{r-1}m^{\hat r-1})$ time complexity.  
\end{corollary}

\vspace{-10pt}

\subsection{Special Case II: \ref{scca} with a rank-one cross-covariance matrix}\label{sub:rankone}
In this subsection, we study the other interesting low-rank special case of  \ref{scca} where the cross-covariance matrix $\bm A$ is rank-one. For this special case, we
prove its NP-hardness with reduction to the sparse regression problem.

We observe that \ref{scca} can be separable over variables $\bm x$ and $ \bm y$ for the rank-one $\bm A$. 
In fact, suppose that $\bm A= \bm a \bm b^{\top}$, then \ref{scca} is equivalent to
\begin{equation} \label{eq_rankone}
v^*:=	\max_{\bm x \in \Re^n, \bm y \in \Re^m} \left\{ \bm x^{\top} \bm a\bm b^{\top} \bm y: \bm x^{\top}\bm B \bm x \le 1, \bm y^{\top}\bm C \bm y \le 1, \|\bm x\|_0 \le s_1, \|\bm y\|_0 \le s_2\right\}
\end{equation}
which  can be equivalently the product of the optimal values of the following two subproblems:
\begin{equation} \label{eq:sep}
\begin{aligned}
v_x:=\max_{\bm x\in \Re^n} \{\bm a^{\top} \bm x:  \bm x^{\top} \bm B \bm x\le 1, \|\bm x\|_0 \le s_1  \},\\
v_y:=\max_{\bm y\in \Re^m} \{\bm b^{\top} \bm y:  \bm y^{\top} \bm C \bm y\le 1, \|\bm y\|_0 \le s_2  \}.
\end{aligned}
\end{equation}
That is, the identity $v^*=v_xv_y$ holds. According to \Cref{prop:bigm},  introducing binary variables, we can reformulate two subproblems \eqref{eq:sdp} as mixed-integer convex quadratic programs. Consequently, the rank-one SCCA problem, as formulated in \eqref{eq_rankone}, simplifies to two mixed-integer convex quadratic programs. This simplification is much more tractable compared to addressing the MISDP \eqref{eq:misdp}, which involves a large-sized positive semidefinite variable $\bm X$ of dimension $(n+m)\times (n+m)$. Our numerical findings confirm the reduced complexity of the rank-one SCCA model. %Thus, the rank-one SCCA \eqref{eq_rankone} is equivalent to solving two mixed-integer convex quadratic programs, which can be easier to address compared to the MISDP \eqref{eq:misdp} involving a sized-$(n+m)\times (n+m)$ positive semidefinite variable $\bm X$. Our numerical study also demonstrates this phenomenon.

Next, we show that each subproblem in  \eqref{eq:sep} can be reduced to the classic sparse regression problem \cite{atamturk2019rank,miller2002subset} and is thus NP-hard as shown below.
\begin{theorem}\label{them:np}
 When matrix $\bm A:=\bm a^{\top}\bm b$ is rank-one, each maximization problem in \eqref{eq:sep} is NP-hard. 
\end{theorem}
 \noindent\textit{Proof.}
    See Appendix A.5. \qed

%\Cref{them:np} connects the maximization problems in \eqref{eq:sep} with the classic sparse regression problem, implying that even the rank-one SCCA \eqref{eq_rankone} is NP-hard to solve.
Theorem \ref{them:np} links the maximization problem \eqref{eq:sep} and the well-known sparse regression problem, implying that even solving the rank-one SCCA problem \eqref{eq_rankone} is NP-hard. This also suggests employing strong perspective formulations (see, e.g., \cite{atamturk2019rank,xie2020scalable}) when solving the subproblems \eqref{eq:sep}, which are shown to be stronger and easier to solve than the SDP relaxation \eqref{eq:sdp} in our numerical study.

\vspace{-10pt}

\section{Numerical results}\label{sec:numerical}
This section tests the numerical performance of our formulations and algorithms on synthetic data.   All the experiments are conducted in Python 3.6 with calls to Gurobi 9.5.2 and MOSEK 10.0.29 on a PC with 10-core CPU, 16-core GPU, and 16GB of memory. 

We generate random instances by fixing the dimensions $n, m$ and the sparsity levels $s_1, s_2$. For each instance, given parameters $(n, m, s_1, s_2)$, %by setting the random seed to be 1, 
we first generate the covariance matrix $\begin{pmatrix}
    \bm B^0 & \bm A^0\\
    (\bm A^0)^{\top} & \bm C^0
\end{pmatrix}$ as follows;
\begin{enumerate}[(i)]
    \item $\bm B^0 \in \S_{++}^n$: Let $\hat{\bm B}$ consist of $n\times n$ elements generated from a normal distribution $\N(0,1)$, and let $\bm B^0=\hat{\bm B} \hat{\bm B}^{\top} + \bm I$;
    \item $\bm C^0 \in \S_{++}^m$: Let $\hat{\bm C}$ consist of $m\times m$ elements generated from a normal distribution $\N(0,1)$, and let $\bm C^0=\hat{\bm C} \hat{\bm C}^{\top} + \bm I$; and
    \item $\bm A^0 \in \Re^{n\times m}:= \lambda \bm B^0 \bm u \bm v^{\top} \bm C^0 $: We generate    $\lambda$ uniformly from $(0, 1)$, and  vectors $\bm u, \bm v$ are generated from a normal distribution $\N(0,1)$ that satisfy $\|\bm u\|_0 = s_1$, $\|\bm v\|_0 = s_2$,   $\bm u^{\top} \bm B^0 \bm u =1 $ and  $\bm v^{\top} \bm C^0 \bm v =1$.
\end{enumerate}
Next, we sample $N=5,000$ data points $\{(\bm u_i, \bm v_i)\}_{i\in [N]} \in \Re^n\times \Re^m$ from the normal distribution with zero mean and the covariance $\begin{pmatrix}
 \bm B^0 & \bm A^0\\
    (\bm A^0)^{\top} & \bm C^0
\end{pmatrix}$.  Then, let us estimate $\bm A^0, \bm B^0, \bm C^0$ by sample covariance matrices below
\[\bm A = \sum_{i\in [N]} \bm u_i\bm v_i^{\top}, \ \ \bm B = \sum_{i\in [N]} \bm u_i\bm u_i^{\top}, \ \ \bm C = \sum_{i\in [N]} \bm v_i\bm v_i^{\top}. \]

The numerical results are presented in Table \ref{table1} that include multiple instances with various parameters $(n,m,s_1, s_2)$. Throughout, the computational time is in seconds, the time limit is one hour, and the dashed line ``-" denotes the unsolved case within the time limit.
First, based on the combinatorial formulation \eqref{scca:com}, we consider using the greedy and local search algorithms to approximately solve \ref{scca}, and their detailed implementation can be found in Appendix B.
Note that we let  \textbf{LB} denote the lower bound obtained from the approximation algorithm. In \Cref{table1}, we define \textbf{gap(\%)}$:=(\hat v - v^*)/v^*$  to be the optimality gap of the upper bound in \eqref{eq:sdp}, and we replace $v^*$ with the best lower bound when  $v^*$ is not available.  It is seen that the greedy and local search algorithms are quite scalable, and the SDP relaxation \eqref{eq:sdp} yields a tight upper bound with an optimality gap at most $8.16\%$. We apply a branch-and-cut algorithm to solve SCCA \eqref{scca:milp} via the delayed cut generation procedure, which can handle the case up to size 20 in \Cref{table1}. One reason may be because  SCCA \eqref{scca:milp} has a weak relaxation bound. Therefore,
although our proposed cut in \Cref{sec:misdp} admits closed form, the branch-and-cut algorithm explores a considerable amount of nodes before termination.

\begin{table}
%\vspace{-15pt}
		\centering
  \scriptsize
  \caption{Solving SCCA with synthetic data}
		\setlength{\tabcolsep}{2pt}\renewcommand{\arraystretch}{1.1}
		\begin{tabular}{c c c c|r r| r  r| r r r | r r}
			\hline 
		 & & & & \multicolumn{2}{c|}{Greedy} & \multicolumn{2}{c|}{Local search} & \multicolumn{3}{c|}{SDP relaxation \eqref{eq:sdp}} & \multicolumn{2}{c}{SCCA \eqref{scca:milp}}
			\\ \cline{5-13}
		$n$ & $m$ & $s_1$ &	 $s_2$&  \multicolumn{1}{c}{LB}  & \multicolumn{1}{c|}{time(s)}   &  \multicolumn{1}{c}{LB}  & \multicolumn{1}{c|}{time(s)}  &\multicolumn{1}{c}{$\hat v$} &\multicolumn{1}{c}{gap(\%)} &\multicolumn{1}{c|}{time(s)} &\multicolumn{1}{c}{$v^*$} &\multicolumn{1}{c}{time(s)} \\
			\hline
   % 5 & 5 & 3 & 3 & 0.8854 &1 & 0.8854 &1 &  0.9210 & 1 &0.9209 & 1\\
   %  10 & 10 & 3 & 3 & 0.4647 &1  & 0.4647 &1 & 0.4672 &1 & 0.4647 &  \\
    10 & 10 & 5 & 5 & 0.244 & 1 & 0.244 & 1& 0.247 & 1.33 &1 & 0.244 & 26\\ 
    20 & 20 & 5 & 5 & 0.244 & 1 & 0.244 & 1& 0.256 & 1.23 & 1 & 0.244 & 2217\\  
20 & 20 & 10 & 10 & 0.275 & 1 & 0.275 & 1& 0.278 & 1.23 & 1 & 0.275 & 3562\\  
40 & 40 & 5 & 5 & 0.695 & 1 & 0.695 & 1 & 0.701 & 0.83 &1 &  -& -\\
40 & 40 & 10 & 10 & 0.705 & 1 & 0.705 & 1 & 0.708 & 0.45 & 1 &  -& -\\
40 & 60 & 5 & 10 & 0.707 & 1 & 0.707 & 1 & 0.714 & 0.93 & 1 &  -& -\\
40 & 60 & 10 & 5 & 0.704 & 1 &  0.704& 1 & 0.708 & 0.65 & 1 &  -& -\\
60 & 60 & 5 & 5 & 0.720& 1 & 0.720 & 1 & 0.727 & 0.86& 14 &  -& -\\
60 & 60 & 10 & 10 &0.714 & 1 & 0.714 & 1 &0.721 & 1.00 &12 &  -& -\\
80 & 80 & 5 & 5 & 0.395& 1 & 0.395& 1 & 0.427 & 8.16 & 56 &  -& -\\
80 & 80 & 10 & 10 & 0.399 & 1 &0.399 & 1 & 0.428& 7.36 & 62 &  -& -\\
100 & 100 & 5 & 5 &0.942 & 1 & 0.942 & 1 & 0.944 & 0.23 & 257 &  -& -\\
100 & 100 & 10 & 10 & 0.940 & 1 & 0.940 & 1 & 0.942 & 0.23 & 313 &  -& -\\
120 & 120 & 5 & 5 & 0.479 & 1 & 0.479 & 1 & 0.517  & 7.90 & 1360 &  -& -\\
120 & 120 & 10 & 10 & 0.501& 1 & 0.501 & 1 & 0.942 & 7.86 & 1569 &  -& -\\
   \hline
		\end{tabular}%
  \label{table1}
  \vspace{-15pt}
\end{table}

The complexity analysis of low-rank \ref{scca} in \Cref{sec:rank} indicates that rank-one SCCA \eqref{eq_rankone} can be more tractable, as we decompose it into two subproblems in \eqref{eq:sep}. By approximating $\bm A$ with a rank-one matrix that consists of leading singular value and vectors, \Cref{table:rankone} presents the numerical results for solving rank-one SCCA \eqref{eq_rankone}. In addition to the SDP relaxation \eqref{eq:sdp}, we consider the strong perspective formulations of subproblems \eqref{eq:sep} to provide an upper bound for rank-one SCCA \eqref{eq_rankone} (see, e.g., \cite{atamturk2019rank,xie2020scalable}), denoted by \textbf{Perspective} in \Cref{table:rankone}. We also compute its optimality gap and compare it with SDP relaxation \eqref{eq:sdp}. It is obvious that perspective relaxation is computationally efficient and yields smaller optimality gaps, which solves all the testing cases in 15 seconds with an optimality gap of up to 11.3\%.
As previously mentioned in \Cref{sub:rankone}, we can solve two mixed-integer quadratic programs below via Gurobi to find the optimal value of rank-one SCCA \eqref{eq_rankone}, i.e., $v^*:=v_xv_y$, where the performance can be found in the last column of \Cref{sub:rankone}. We see that we can solve size-$100\times 100$ rank-one SCCA \eqref{eq_rankone}.
% \begin{equation} \label{eq:miqp}
% \begin{aligned}
% v_x:=\max_{\bm x\in \Re^n} \{\bm a^{\top} \bm x:  \bm x^{\top} \bm B \bm x\le 1, \|\bm x\|_0 \le s_1  \},\\
% v_y:=\max_{\bm y\in \Re^m} \{\bm b^{\top} \bm y:  \bm y^{\top} \bm C \bm y\le 1, \|\bm y\|_0 \le s_2  \}.
% \end{aligned}
% \end{equation}
\begin{table}[htb]
		\centering  \scriptsize
  \caption{Solving rank-one SCCA with synthetic data}
		\setlength{\tabcolsep}{1.5pt}\renewcommand{\arraystretch}{1.2}
		\begin{tabular}{c c c c|r r| r r| r  r | r r | r r}
			\hline 
		 & & & & \multicolumn{2}{c|}{Greedy} & \multicolumn{2}{c|}{Local Search} & \multicolumn{2}{c|}{SDP relaxtion \eqref{eq:sdp}} & \multicolumn{2}{c|}{Perspective} & \multicolumn{2}{c}{SCCA \eqref{eq_rankone}}
			\\ \cline{5-14}
		$n$ & $m$ & $s_1$ &	 $s_2$&  \multicolumn{1}{c}{LB}  & \multicolumn{1}{c|}{time(s)}    &  \multicolumn{1}{c}{LB}  & \multicolumn{1}{c|}{time(s)}  &\multicolumn{1}{c}{gap(\%)} &\multicolumn{1}{c|}{time(s)} &\multicolumn{1}{c}{gap(\%)} &\multicolumn{1}{c|}{time(s)}  &\multicolumn{1}{c}{$v^*$} &\multicolumn{1}{c}{time(s)}   \\
			\hline
   % 5 & 5 & 3 & 3 & 0.8777 &1 & 0.8777 &1 & 0.9095 & 1 &0.9093 & 1 & 0.9091 & 1\\
   %  10 & 10 & 3 & 3 &  0.4745 &1  &  0.4745 &1 & 0.4748 &1 & 0.4745 & 1 &0.4745 & 1\\
%     10 & 10 & 4 & 4 & 0.4717 & 1 & 0.4717 & 1& 0.4719 & 1 & 0.4717 & 1 & 0.4717 & 1 \\ 
% 15 & 15 & 5 & 5 & 0.4114 & 1 & 0.4114 & 1& 0.4127 & 1 & 0.4127 & 1 & 0.4121 & 1\\ 
%    30 &  15 & 10 & 5 & 0.9071 & 1 & 0.9071 & 1 & 0.9078 & 1 \\
%    45 & 30 & 10 & 10 & 0.3578 & 1 & 0.3578 & 1 & 0.3685& 1\\
% 35 & 35 & 10 & 10 & 0.9711 & 1 & 0.9716 & 1 & 0.9720 & 1\\
% 45 & 45 & 5 & 5 & 0.9307 & 1 & 0.9307 & 1 & 0.9319 & 3 \\
% 45 & 45 & 10 & 10 & 0.9305 & 1 & 0.9305 & 1 & 0.9316& 3 \\
% 50 & 50 & 5 & 10 & 0.5771 & 1 & 0.5771 & 1 & 0.5818& 8 & 0.5816 & 1 & 0.5771 & 4\\

% 50 & 50 & 10 & 5 & 0.5676 & 1 & 0.5676 & 1 & 0.5726 & 7 & 0.5723 & 1 & 0.5676 & 10 \\
% 50 & 50 & 10 & 10 & 0.5682 & 1 & 0.5682 & 1 & 0.5727 & 7 & 0.5725 & 1 & 0.5683 & 60\\
% 80 & 80 & 5 & 10 & 0.5771 & 1 & 0.5771 & 1 & 0.5818& 8 & 0.5816 & 1 & 0.5771 & 4\\

50 & 50 & 10 & 10 & 0.382 & 1 & 0.382 & 1 & 3.79 & 6 & 2.44 & 1 & 0.382 & 30 \\
50 & 50 & 20 & 20 & 0.409 & 1 &0.409 & 1 & 2.81 & 7 & 1.74 & 1 & 0.409 & 293\\
100 & 100 & 10 & 10 & 0.928 & 1 & 0.928 & 1 & 0.79  & 492 &  0.47 & 1 & 0.928 & 81\\
100 & 100 & 20 & 20 & 0.943 & 1 & 0.943 & 2 &  0.49 & 685 & 0.31 & 1 & 0.943 & 3463\\
200 & 200 & 10 & 10 & 0.549 & 1 & 0.549 & 1 &  -& -&  7.38 & 1 & -& - \\
200 & 200 & 20 & 20 & 0.524 & 1 & 0.524& 5 &  -& - & 9.70 & 1 & -& -\\
300 & 300 & 10 & 10 & 0.874 & 1 & 0.874 & 1 &  -& - &  2.56 & 6 & -& -\\
300 & 300 & 20 & 20 & 0.878 & 1 & 0.878 & 9 &  -&- & 2.49 & 8 & -& -\\
400 & 400 & 10 & 10 & 0.840 & 1 & 0.840 & 2 &  -& -&  4.43 & 9 & -& -\\
400 & 400 & 20 & 20 & 0.842 & 1 & 0.842 & 14 &  -& - & 4.34 & 10 & -& -\\
500 & 500 & 10 & 10 &  0.701 & 1 & 0.701 & 2 &  - & - &  11.3& 14 & -& -\\
500 & 500 & 20 & 20 & 0.710 & 6 & 0.710  & 59 &  - & - & 10.9 & 15 & -& -\\
% 600 & 600 & 10 & 10 & 0.9286 & 1 & 0.9286 & 1 &  0.9491 & 650 &  0.9481 & 1\\
% 600 & 600 & 20 & 20 & 0.9289 & 1 & 0.9289 & 1 &  0.9450 & 836 & 0.9443 & 1\\
% 700 & 700 & 10 & 10 & 0.9286 & 1 & 0.9286 & 1 &  0.9491 & 650 &  0.9481 & 1\\
% 700 & 700 & 20 & 20 & 0.9289 & 1 & 0.9289 & 1 &  0.9450 & 836 & 0.9443 & 1\\
% 900 & 900 & 10 & 10 & 0.92 & 1 & 0.9286 & 1 &  -& 3600 &  0.94 & 1\\
% 900 & 900 & 20 & 20 & 0.928 & 1 & 0.9289 & 1 &  -& 3600 & 0.94 & 1\\
   \hline
		\end{tabular}%
  \label{table:rankone}
\end{table}

\noindent\textbf{Acknowledgements:} The authors would like to thank Rahul Mazumder for introducing the problem to us and for sharing data and references.

{\bibliography{reference.bib}}

\newpage
\section*{Appendix A: Proofs}
\subsection*{A.1 \ \ Proof of \Cref{prop:cca}}
 \noindent\textit{Proof.}
The proof includes three parts.

\noindent\textbf{Part (i).} To prove the equivalence between \ref{cca} and its \ref{sdp}, let us introduce the Lagrangian multiplies $\theta_1\geq 0,\theta_2 \ge 0$ corresponding to two constraints in \ref{sdp}, which leads to the following
Lagrangian dual problem
\begin{align}\label{eq:min}
&\min_{\theta_1\ge 0, \theta_2\ge 0}  \left\{\theta_1 + \theta_2: 
   \theta_1 \tilde{\bm B} +  \theta_2 \tilde{\bm C} \succeq   \tilde{\bm A} \right\}=
\min_{\theta_1\ge 0, \theta_2\ge 0}  \left\{\theta_1 + \theta_2: \begin{pmatrix}
   \theta_1 \bm B & \frac{\bm A}{-2}\\
\frac{\bm A^{\top}}{-2} & \theta_2 \bm C
\end{pmatrix} \succeq 0\right\} 
\end{align}
where the  equation results from the definition of block matrices $\tilde{\bm A}, \tilde{\bm B},$ and $\tilde{\bm C}$. Given the nonzero matrices $\bm A\neq \bm 0, \bm B\neq \bm 0, \bm C \neq \bm 0$ and positive semidefinite matrices $\bm B \succeq 0, \bm C\succeq 0$, following \Cref{lem:psd}, we must have $\theta_2\bm C - \bm A^{\top}(\theta_1\bm B)^{\dag} \bm A /4\succeq 0$ and $\theta_1\bm B - \bm A(\theta_2\bm C)^{\dag} \bm A^{\top} /4\succeq 0$, implying that either $\theta_1=0$ or $\theta_2=0$ is infeasible to the minimization problem above. That is, $\theta_1>0$ and $\theta_2>0$ must hold. 

According to  \Cref{lem:psd}, the block matrix $\begin{pmatrix}
    \bm B & \bm A\\
    \bm A^{\top} & \bm C
\end{pmatrix}$  is positive semidefinite, implying that $(\bm I-\bm C \bm C^{\dag})\bm A ^{\top}=\bm 0, (\bm I-\bm B \bm B^{\dag})\bm A =\bm 0$. Then, it is easy to show 
$$\left(\bm I-\theta_2\bm C (\theta_2\bm C)^{\dag}\right)\frac{\bm A^{\top}}{2}=\bm 0, \forall \theta_2>0.$$ 
Given $\theta_1, \theta_2>0$ and using \Cref{lem:psd}, the result above allows  us to further simplify the right-hand side minimization problem in \eqref{eq:min} to
\begin{align*}
&\min_{\theta_1\ge 0, \theta_2\ge 0}  \left\{\theta_1 + \theta_2:  4\theta_1\theta_2  \bm B \succeq  \bm A \bm C^{\dag} \bm A^{\top}  \right\} \\
&= \min_{\theta_1\ge 0, \theta_2\ge 0}  \left\{\theta_1 + \theta_2:  4\theta_1\theta_2 \ge \sigma^2_{\max}\left(\sqrt{\bm B^{\dag}} \bm A\sqrt{\bm C^{\dag}}\right) \right\}  = \sigma_{\max}\left(\sqrt{\bm B^{\dag}} \bm A \sqrt{\bm C^{\dag}}\right),
\end{align*}
where  the first equation is because
\begin{align*}
  &4\theta_1\theta_2  \bm B \succeq  \bm A \bm C^{\dag} \bm A^{\top} \Longleftrightarrow
4\theta_1\theta_2  \bm I \succeq \sqrt{\bm \Lambda^{-1}} \bm Q^{\top} \bm A \bm C^{\dag} \bm A^{\top}  \bm Q \sqrt{\bm \Lambda^{-1}}\\ &\Longleftrightarrow  4\theta_1\theta_2 \ge \lambda_{\max}\left( \sqrt{\bm \Lambda^{-1}} \bm Q^{\top} \bm A \bm C^{\dag} \bm A^{\top}  \bm Q \sqrt{\bm \Lambda^{-1}} \right) \\
  & \Longleftrightarrow  4\theta_1\theta_2 \ge \lambda_{\max}\left( \sqrt{{\bm C}^{\dag}}   \bm A^{\top}  \bm B^{\dag} \bm A  \sqrt{\bm C^{\dag}} \right)  \Longleftrightarrow 4\theta_1\theta_2 \ge \sigma^2_{\max}\left(\sqrt{\bm B^{\dag}} \bm A\sqrt{\bm C^{\dag}}\right),
\end{align*}
where we let $\bm B = \bm Q\bm \Lambda \bm Q^{\top}$ denote the eigendecomposition of matrix $\bm B$ with $\bm \Lambda$ containing all the positive eigenvalues.

As a result, the dual problem of \ref{sdp}  admits an optimal value of $\sigma_{\max}\left(\sqrt{\bm B^{\dag}} \bm A \sqrt{\bm C^{\dag}}\right)$, which gives an upper bound of the \ref{cca} and its \ref{sdp}. Next, we construct their optimal solutions, which exactly attain this upper bound. Thus, this upper bound 
% $\sigma_{\max}\left(\sqrt{\bm B^{\dag}} \bm A \sqrt{\bm C^{\dag}}\right)$ 
is achievable and equals their optimal values.
 
\noindent\textbf{Part (ii).} For the \ref{cca},
let us consider a part of optimal solutions $(\bm x^*, \bm y^*)$ below
 $$\bm x^* = \sqrt{\bm B^{\dag}} \bm q, \ \  \bm y^* = \sqrt{\bm C^{\dag}} \bm p,$$ 
  with $\bm q \in \Re^n, \bm p \in \Re^m$ denoting a pair of leading singular vectors of matrix $\sqrt{\bm B^{\dag}} \bm A \sqrt{\bm C^{\dag}}$.

    First, $(\bm x^*, \bm y^*)$ is feasible to the \ref{cca} as
    \[(\bm x^*)^{\top} \bm B \bm x^*= \bm  q^{\top} \sqrt{\bm B^{\dag}} \bm B \sqrt{\bm B^{\dag}} \bm q \le \bm q^{\top} \bm q = 1, \ \ (\bm y^*)^{\top} \bm C \bm y^*= \bm  p^{\top} \sqrt{\bm C^{\dag}} \bm C \sqrt{\bm C^{\dag}} \bm p \le \bm p^{\top} \bm p = 1, \]
    where the inequalities stem from the facts that $ \bm I \succeq \sqrt{\bm B^{\dag}} \bm B \sqrt{\bm B^{\dag}} $ and $ \bm I \succeq \sqrt{\bm C^{\dag}} \bm C \sqrt{\bm C^{\dag}} $.
    
    On the other hand, according to the definitions of $\bm q, \bm p$, we can show that $(\bm x^*, \bm y^*)$ is optimal to the \ref{cca}, i.e.,
    $$(\bm x^*)^{\top} \bm A \bm y^*= \bm  q^{\top} \sqrt{\bm B^{\dag}} \bm A \sqrt{\bm C^{\dag}} \bm p=  \sigma_{\max}\left(\sqrt{\bm B^{\dag}} \bm A \sqrt{\bm C^{\dag}} \right).$$ 
\noindent\textbf{Part (iii).} In a similar vein, we can show that $\bm X^* =  \begin{pmatrix}
    \bm x^*\\
    \bm y^*
\end{pmatrix} \begin{pmatrix}
    \bm x^*\\
    \bm y^*
\end{pmatrix}^{\top}$ is  optimal to \ref{sdp} with the optimal value $\sigma_{\max}\left(\sqrt{\bm B^{\dag}} \bm A \sqrt{\bm C^{\dag}} \right)$. 
\qed

%{Different from previous studies on CCA. By exploiting the structure}

% \subsection*{A.2 \ \ Proof of \Cref{them:com}}
%  \noindent\textit{Proof.}
% By introducing the  subsets $(S_1, S_2)$ to denote the supports of variables $(\bm x, \bm y)$ in \ref{scca},  then we can remove the zero-norm constraints on $(\bm x, \bm y)$ and reformulate \ref{scca} as
%  \begin{align} \label{scca_set}
%      &v^*:= \max_{\begin{subarray}{c}
% 			S_1\subseteq [m], |S_1|\le s_1, \\
% 			S_2 \subseteq [n], |S_2|\le s_2
% 			\end{subarray}} \max_{
%    \begin{subarray}{c}
%    \bm x \in \Re^{|S_1|},\\ \bm y \in \Re^{|S_2|}
%    \end{subarray}} \left\{ \bm x^{\top} \bm A_{S_1, S_2} \bm y:  \bm x^{\top}\bm B_{S_1, S_1} \bm x \le 1, \bm y^{\top}\bm C_{S_2, S_2} \bm y \le 1\right\}. 
%  \end{align}

% Following from the Part (i) in \Cref{prop:cca}, we can show that
%  for any subsets $S_1\subseteq [n], S_2\subseteq [m]$, the following identity holds.
%       \begin{align*}
%     & \max_{   \begin{subarray}{c}
%    \bm x \in \Re^{|S_1|}, \bm y \in \Re^{|S_2|}
%    \end{subarray}} \left\{ \bm x^{\top} \bm A_{S_1, S_2} \bm y:  \bm x^{\top}\bm B_{S_1, S_1} \bm x \le 1, \bm y^{\top}\bm C_{S_2, S_2} \bm y \le 1\right\} \\
%    &= \sigma_{\max}\left(\sqrt{(\bm B_{S_1, S_1})^{\dag}} \bm A_{S_1, S_2} \sqrt{(\bm C_{S_2, S_2})^{\dag}}\right).
%  \end{align*}
% Plugging the result above into the inner maximization problem in \eqref{scca_set}, we complete the proof.
% \qed

\subsection*{A.2 \ \ Proof of \Cref{prop:bigm}}
 \noindent\textit{Proof.}
Let $(\bm x^*, \bm y^*)$ denote an optimal solution to \ref{scca}. We bound $\|\bm x^*\|_2$ first and the same technique can be also straightforwardly applied to bound $\|\bm y^*\|_2$.

For matrix $\bm B\in \S_+^n$ of rank $r$, we let $\{\bm q_i\}_{i\in [n]} \in \Re^n$ denote the  eigenvectors corresponding to $n$ eigenvalues $\bm\lambda$ of $\bm B$ such that $\lambda_1\geq \ldots \geq\lambda_r> \lambda_{r+1}=\ldots=\lambda_{n}=0$. Thus, $\{\bm q_i\}_{i\in [n]}$ are orthonormal and span the space of $\Re^n$. Hence, there exists $\bm\alpha\in \Re^n$ such that
$\bm x^*=\sum_{i\in [n]}\alpha_i\bm q_i$.
Given that $(\bm x^*)^\top \bm B \bm x^*\leq 1$, we have
\begin{align*}
\sum_{i\in [r]}\alpha_i^2\lambda_i\leq 1.
\end{align*}
Hence, the values of $\{\alpha_i\}_{i\in [r]}$ are bounded. On the other hand, let us define a subset $S\subseteq [n]$ of size at most $s_1$ such that $x_i^*\neq 0$ for each $i\in S$ and $x_j^*=0$ for each $j\in [n]\setminus S$. Then for each $j\in [n]\setminus S$, we arrive at the following linear system:
\begin{align}
\sum_{j\in [r+1,n]}\alpha_j\hat{\bm q}_j=-\sum_{i\in [r]}\alpha_i\hat{\bm q}_i,\label{eq_alpha_i}
\end{align}
where $\hat{\bm q}_i$ denote a subvector of ${\bm q}_i$ with indices $[n]\setminus S$ for each $i\in [n]$. For a fixed $\{\alpha_i\}_{i\in [r]}$, since the linear system \eqref{eq_alpha_i} is nonempty, we let 
$\bar{\bm Q} \bar{\bm \alpha}=\bar{\bm q}$ denote  its minimal linear subsystem such that a submatrix $\bar{\bm Q}$ is non-singular and the index set $\hat{S}$ of $\bar{\bm \alpha}$ is a subset of $[n]\setminus S$. Thus, we can construct an alternative solution $\hat{\bm\alpha}$ such that
\[\hat{\alpha}_i=\begin{cases}
\alpha_i, &\text{ if }i\in [r],\\
(\bar{\bm Q} ^{-1}\bar{\bm q})_i, &\text{ if }i\in \hat{S},\\
0, &\text{ otherwise},
\end{cases}\]
and $\hat{\bm x}=\sum_{i\in [n]}\hat{\alpha}_i\bm q_i$. According to \Cref{lem:psd}, we have
\begin{align*}
\hat{\bm x}^\top \bm{B} \hat{\bm x}\leq 1,  \hat{\bm x}^\top \bm{A}\bm{y}^*=
({\bm x}^*)^\top \bm{A}\bm{y}^*,
\end{align*}
i.e., $(\hat{\bm x},\bm{y}^*)$ is also optimal to \ref{scca}. Hence,
\begin{align*}
\|\hat{\bm x}\|_2\leq \sqrt{\|\bar{\bm Q} ^{-1}\bar{\bm q}\|_2^2+\sum_{i\in[r]}\alpha_i^2}
\end{align*}
Note that $\sum_{i\in[r]}\alpha_i^2\leq 1/\lambda_r$ and
\begin{align*}
\|\bar{\bm Q} ^{-1}\bar{\bm q}\|_2^2\leq \|\bar{\bm Q} ^{-1}\|_2^2\|\bar{\bm q}\|_2^2\leq 
\frac{1}{s_{\min}(\bm{B})}\frac{1}{\lambda_r}
\end{align*}
where $s_{\min}(\bm{B})$ denotes the smallest nonzero singular values of all the submatrices of $[\bm{q}_{r+1},\ldots,\bm{q}_{n}]$. In summary, we have
\begin{align*}
\|\hat{\bm x}\|_2\leq \sqrt{1/\lambda_r+1/(\lambda_rs_{\min}(\bm{B}))}.
\end{align*}
This completes the proof.
\qed

\subsection*{A.3 \ \ Proof of \Cref{prop:milp}}
 \noindent\textit{Proof.}
First, for any binary variable ${\bm z} \in \Z$, suppose  ${S}_1:=\{i: z_i=1,\forall i\in [n]\}$, ${S}_2:=\{i-n: z_i=1,\forall i\in [n+1,n+m]\}$, and $T \subseteq [n+m]$ denotes the support of $\bm z$. Then following the proof of \Cref{prop:cca}, we can construct a rank-one optimal solution $\bm X^*:=\begin{pmatrix}
    \bm x^*\\
    \bm y^*
\end{pmatrix} \begin{pmatrix}
    \bm x^*\\
    \bm y^*
\end{pmatrix}^{\top}
$ to the maximization problem  below that admits the optimal value $\sigma_{\max}\left(\sqrt{(\bm B_{S_1, S_1})^{\dag}} \bm A_{S_1, S_2} \sqrt{(\bm C_{S_2, S_2})^{\dag}}\right)$, i.e.,
\begin{align*}
&\max_{\bm X \in \S_+^{n+m}} \{\tr(\tilde{\bm A} \bm X): \tr( \tilde{\bm B} \bm X) \le 1, 
 \tr(\tilde{\bm C} \bm X)\le 1,  X_{ii}=0,  \forall i\in [n+m]\setminus T\}\\
 &= \sigma_{\max}\left(\sqrt{(\bm B_{S_1, S_1})^{\dag}} \bm A_{S_1, S_2} \sqrt{(\bm C_{S_2, S_2})^{\dag}}\right) \ge 	f(\bm z),
\end{align*}
where the inequality is because the maximization problem above relaxes the valid constraints $X_{ii} \le M_{ii}$ for all $i\in T$   in maximization problem \eqref{eq:obj}. The result in \Cref{cor_prop:bigm} suggests that $\bm x^*, \bm y^*$ can be bounded and their two norms must not exceed $M_1, M_2$, which means that the optimal solution $\bm X^*$ satisfies the  $X_{ii} \le M_{ii}$ for all $i\in T$. Therefore,  $\bm X^*$ is feasible and optimal to maximization problem \eqref{eq:obj} and we have that  
\begin{align*}
 f(\bm z):= \sigma_{\max}\left(\sqrt{(\bm B_{S_1, S_1})^{\dag}} \bm A_{S_1, S_2} \sqrt{(\bm C_{S_2, S_2})^{\dag}}\right).
\end{align*}

According to strong duality, the minimization problem \eqref{eq:dual} admits an optimal value $\sigma_{\max}\left(\sqrt{(\bm B_{S_1, S_1})^{\dag}} \bm A_{S_1, S_2} \sqrt{(\bm C_{S_2, S_2})^{\dag}}\right)$. Next, we construct its optimal solution $(\theta_1^*, \theta_2^*, \bm \lambda^*)$.
% \begin{align}\label{eq:dual1}
%    f(\bm z):=   \min_{\begin{subarray}{c}
% 	     \theta_1 \ge 0, \theta_2 \ge 0, \\ \bm \lambda \in \Re_+^{n+m}
% 	\end{subarray}}	\bigg\{ \theta_1 + \theta_2 + \sum_{i\in [n+m]} \lambda_i M_{ii}  z_i: \begin{pmatrix}
% 	   & \theta_1\bm B  & -\bm A/2\\
%      & -\bm A^{\top}/2 & \theta_2\bm C
% 	\end{pmatrix} \succeq - \Diag(\bm \lambda) \bigg\},
% \end{align}

For any given $\epsilon>0$, we let $\theta^*_1=f(\bm z)/2$, $\theta^*_2=f(\bm z)/2$, $\hat{\lambda}_i(\epsilon) = \frac{\epsilon}{M_{ii}|T|}$ for all $i\in T$, and $\hat{\lambda}_i(\epsilon) = \lambda^*(\epsilon)$ for all $i\in [n]\setminus T$, where
$$ \lambda^*(\epsilon) := \left[\lambda_{\max} \left( \bm D_2^{\top} \left(\bm D_1 + \Diag\left(\hat{\bm \lambda}_{T}(\epsilon)\right)\right )^{-1} \bm D_2-\bm D_3 \right)\right]_+.$$
 It is easy to compute that $ \theta^*_1 + \theta^*_2 + \sum_{i\in [n+m]} \hat{\lambda}_i(\epsilon) M_{ii}  z_i = f(\bm z)+ \epsilon$.
Thus, for any $\epsilon>0$,  if $(\theta_1^*, \theta_2^*, \hat{\bm \lambda}(\epsilon))$ were feasible, then it is an $\epsilon$-optimal solution to the minimization problem \eqref{eq:dual}. 
It remains to verify the feasibility of the solution  $(\theta_1^*, \theta_2^*, \hat{\bm \lambda}(\epsilon))$, i.e., checking the constraint below
\begin{align*}
    \begin{pmatrix}
	   & \theta_1^*\bm B  & -\bm A/2\\
     & -\bm A^{\top}/2 & \theta_2^*\bm C
	\end{pmatrix} + \Diag\left(\hat{\bm \lambda}(\epsilon)\right)\succeq 0.
\end{align*}
By performing the permutation of the rows and columns of the above matrix, it is sufficient to show that the new block matrix
\begin{align}\label{eq:constr}
\begin{pmatrix}
\bm D_1 + \Diag\left(\hat{\bm \lambda}_{T}(\epsilon)\right) & \bm D_2\\
\bm D_2^{\top} & \bm D_3 + {\lambda^*}(\epsilon) \bm I
\end{pmatrix} \succeq 0,
\end{align}
is positive semidefinite.

Since $\begin{pmatrix}
	   & \bm B_{S_1, S_1}  & -\bm A_{S_1, S_2}/2\\
     & -\bm A_{S_1, S_2}^{\top}/2 & \bm C_{S_2, S_2}
	\end{pmatrix}$ is a principal submatrix of a positive semidefinite matrix $\begin{pmatrix}
	   & \bm B  & -\bm A/2\\
     & -\bm A^{\top}/2 & \bm C
	\end{pmatrix}$, it is also positive semidefinite. According to \Cref{lem:psd} and the fact that $\theta^*_1 = \theta^*_2=\sigma_{\max}\left(\sqrt{(\bm B_{S_1, S_1})^{\dag}} \bm A_{S_1, S_2} \sqrt{(\bm C_{S_2, S_2})^{\dag}}\right)/2$, the matrix $\bm D_1$ is also positive semidefinite. As $\epsilon>0$, the matrix $\bm D_1 + \Diag\left(\hat{\bm \lambda}_{T}(\epsilon)\right)$ must be positive definite, which means that
\[ \left(\bm I - \left(\bm D_1 + \Diag\left(\hat{\bm \lambda}_{T}(\epsilon)\right)\right )\left(\bm D_1 + \Diag\left(\hat{\bm \lambda}_{T}(\epsilon)\right)\right )^{-1}\right )\bm D_2 = \bm 0.\]
Besides, according to the definition of $\lambda^*(\epsilon)$, we obtain 
\[\bm D_3 + \lambda^*(\epsilon) \bm I - \bm D_2^{\top} \left(\bm D_1 + \Diag\left(\hat{\bm \lambda}_{T}(\epsilon)\right)\right )^{-1} \bm D_2\succeq \bm 0.\]
Taking these results together, according to \Cref{lem:psd}, the constraint in \eqref{eq:constr} must hold for a given solution $(\theta_1^*, \theta_2^*, \hat{\bm \lambda}(\epsilon))$. Since the objective value corresponding to $(\theta_1^*, \theta_2^*, \hat{\bm \lambda}(\epsilon))$ is at most $\epsilon$ larger than the optimal value of   problem \eqref{eq:dual}, letting  $\epsilon\to 0$ and using the closedness of the feasible set in problem \eqref{eq:dual}, we can confirm the optimality of $(\theta_1^*, \theta_2^*, \bm \lambda^*)$ with $\lambda_i^*=0$ for all $i\in T$ and $\lambda_i^* = \lambda^*$ for all $i\in [n]\setminus T$.

Given the closed-form optimal solution to problem \eqref{eq:dual}, the rest of the proof follows from \cite{li2021beyond}[theorem 7].
\qed

\subsection*{A.4 \ \ Proof of \Cref{them:rank}}
 \noindent\textit{Proof.}
The proof is split into three parts. \par 
\textbf{Part (i).} It suffices to prove that \ref{cca} admits an optimal solution $(\bm x^*, \bm y^*)$ satisfying  $\|\bm x^*\|_0\le r$ and  $\|\bm y^*\|_0\le \hat r$. Then, $(\bm x^*, \bm y^*)$ is also feasible and optimal to \ref{scca}, which implies the equivalence between \ref{scca} and \ref{cca}.

First, according to Part (ii) in \Cref{prop:cca}, we can obtain a closed-form optimal solution $(\hat{\bm x}, \hat{\bm y})$ for the \ref{cca}. By adjusting $(\hat{\bm x}, \hat{\bm y})$, we will construct a new optimal sparse solution $(\bm x^*, \bm y^*)$ satisfying  $\|\bm x^*\|_0\le r$ and  $\|\bm y^*\|_0\le \hat r$.
     
For matrix $\bm B\in \S_+^n$, we let $\{\bm q_i\}_{i\in [n-r]} \in \Re^n$ denote the  eigenvectors corresponding to $(n-r)$ zero eigenvalues of $\bm B$. Thus, $\{\bm q_i\}_{i\in [n-r]}$ are orthonormal. 
There exists a size-$(n-r)$ subset $S\subseteq [n]$ such that the subvectors $\{(\bm q_i)_S\}_{i\in [n-r]}$ are linearly independent, where $(\bm q_i)_S$ denotes the subvector of $\bm q_i$ indexed by $S$ for each $i\in [n-r]$. %Then, the following $(n-r+1)$ subvectors in $\Re^{n-r}$ must be  linearly dependent
%\[ \hat{\bm x}_S, (\bm q_{1})_S, \cdots, (\bm q_{n-r})_S.\]
As a result, there exist a vector $(\gamma_1, \cdots, \gamma_{n-r})^\top $ such that
\begin{equation}\label{eq:zero}
    \hat{\bm x}_S = \sum_{i\in [n]} \gamma_i (\bm q_{i})_S.
\end{equation}
%where $\alpha\neq 0$ holds because   $\{(\bm q_i)_S\}_{i\in [n-r]}$ are linearly independent. 

%By defining $\beta_i = \gamma_i/\alpha$ for each $i\in [n-r]$, 
Let us now construct solution $\bm x^*$
\[\bm x^* = \hat{\bm x} - \sum_{i\in [n-r]} \gamma_i \bm q_i,\]
where $x^*_i=0$ for all $i\in S$ based on the equation \eqref{eq:zero} and $|S|=n-r$, implying $\|\bm x^*\|_0 \le r$. In addition, we show that the new solution $\bm x^*$ is still optimal to \ref{cca}. First, $\bm x^*$ is feasible since
\[(\bm x^*)^{\top} \bm B (\bm x^*) = \hat{\bm x}^{\top} \bm B  \hat{\bm x}\le 1, \]
where the equation is due to $\bm B \bm q_i = \bm 0$ for all $i\in [n-r]$.

Given the positive semidefinite block matrix $\begin{pmatrix}
    \bm B & \bm A\\
    \bm A^{\top} & \bm C
\end{pmatrix}$, 
using Part (ii) of \Cref{lem:psd}, the identity $(\bm I - \bm B \bm B^{\dag}) \bm A = \bm 0$ is equivalent to $\sum_{i\in [n-r]} \bm q_i  \bm q_i^{\top} \bm A = \bm 0$. Then, for each $i\in [n-r]$, multiplying $\bm q_i^{\top}$ on both sides of this equation leads to
$$\bm q_i^{\top} \bigg(\sum_{j\in [n-r]} \bm q_j  \bm q_j^{\top} \bm A\bigg) \bm A =  \bm q_i^{\top} \bm 0  \Longrightarrow \bm q_i^{\top} \bm A = \bm 0,$$
where the result follows from $\bm q_i^{\top} \bm q_j=0$ for any $i\neq j$.
Then, we can show the optimality of the new solution $\bm x^*$:
\[(\bm x^*)^{\top} \bm A \hat{\bm y} = \hat{\bm x}^{\top} \bm A  \hat{\bm y} + \sum_{i\in [n-r]} \beta_i \bm q_i^{\top}  \bm A \hat{\bm y}= \hat{\bm x}^{\top} \bm A \hat{\bm y}. \]

Similarly, we can also construct an optimal sparse solution $\bm y^*$ by leveraging $\hat{\bm y}$ and eigenvectors of zero eigenvalues of $\bm C$ such that $\|\bm y^*\|_0 \le s_2$.

Therefore, there exists an optimal solution $(\bm x^*, \bm y^*)$ to the \ref{cca} whose zero norms are bounded from above by $r, \hat r$, respectively. 
Adding the constraints $\|\bm x\|_0\le r, \|\bm y\|_0\le \hat r$ to the \ref{cca} does not affect the optimality, which gives an equivalent formulation \eqref{scca_re} of \ref{cca}.  

\textbf{Part (ii).} Suppose that $(\tilde{\bm x}, \tilde{\bm y})$ denotes an optimal solution to problem \eqref{scca_re2}.
When $s_1 \ge r$, following the proof of Part (I), $\tilde{\bm x}$, we can construct another optimal solution $\bm{x}^*$ whose zero norm is bounded by $r$ and $(\bm{x}^*, \tilde{\bm y})$ is feasible and optimal to \ref{scca}.
% \begin{align*}
%     \max_{\bm x \in \Re^n, \bm y \in \Re^m} \left\{ \bm x^{\top} \bm A \bm y: \bm x^{\top}\bm B \bm x \le 1, \bm y^{\top}\bm C \bm y \le 1, \|\bm x\|_0 \le \min\{s_1, r\}, \|\bm y\|_0 \le s_2\right\}.
% \end{align*}

\textbf{Part (iii).} Similarly, we can reduce \ref{scca} to problem \eqref{scca_re3}. We thus complete the proof. \qed

\subsection*{A.5 \ \ Proof of \Cref{them:np}}
 \noindent\textit{Proof.}
Let us first consider the maximization problem over $\bm x$ in \eqref{eq:sep}, i.e.,
\begin{equation} \label{eq:sepx}
\begin{aligned}
v_x:=\max_{\bm x\in \Re^n} \{\bm a^{\top} \bm x:  \bm x^{\top} \bm B \bm x\le 1, \|\bm x\|_0 \le s_1  \}.
\end{aligned}
\end{equation}
Then, we derive a combinatorial optimization reformulation of problem \eqref{eq:sepx} based on the result below.
\begin{claim1}\label{claim1}
For any subset $S\subseteq [n]$, $\max_{\bm x\in \Re^{|S|}} \{\bm a_{S}^{\top} \bm x: \bm x^{\top} \bm B_{S, S} \bm x\le 1 \} = \sqrt{\bm a_{S}^{\top} (\bm B_{S, S})^{\dag} \bm a_{S}}$.
\end{claim1}
 \noindent\textit{Proof.}
Given $\bm A=\bm a\bm b^{\top}$, since the matrix $\begin{pmatrix}
    \bm B & \bm a\bm b^{\top}\\
    \bm b^{\top}\bm a & \bm C
\end{pmatrix}$ is positive semidefinite, using \Cref{lem:psd}, the identity $(\bm I - \bm B_{S,S}\bm B_{S,S}^{\dag})\bm a_S\bm b^{\top} = \bm 0$ must hold for any subset $S$. As a result, we have $\bm a_S - \bm B_{S,S}\bm B_{S,S}^{\dag}\bm a_S= \bm 0$ as vector $\bm b$ is nonzero.

Next, the Lagrangian dual of the problem $\max_{\bm x\in \Re^{|S|}} \{\bm a_{S}^{\top} \bm x: \bm x^{\top} \bm B_{S, S} \bm x\le 1 \}$ can be written as 
\begin{align*}
    \max_{\bm x\in \Re^{|S|}} \{\bm a_{S}^{\top} \bm x: \bm x^{\top} \bm B_{S, S} \bm x\le 1 \} &= \min_{\mu\ge 0} \max_{\bm x\in \Re^{|S|}} \bm a_{S}^{\top} \bm x + \mu-\mu  \bm x^{\top} \bm B_{S, S} \bm x\\
    &= \min_{\mu\ge 0} \mu + \frac{\bm a_S^{\top}\bm B^{\dag}_{S,S}\bm a_S}{4\mu} =  \sqrt{\bm a_{S}^{\top} (\bm B_{S, S})^{\dag} \bm a_{S}},
\end{align*}
where the second equation builds on the identity $\bm a_S - \bm B_{S,S}\bm B_{S,S}^{\dag}\bm a_S= \bm 0$ and optimal solution $\bm x^* = \frac{\bm B^{\dag}_{S,S}\bm a_S}{\sqrt{\bm a_{S}^{\top} (\bm B_{S, S})^{\dag} \bm a_{S}}}$.
    \qedA

Suppose that an optimal solution to problem  \eqref{eq:sepx} admits the support $S^*$. According to Claim \ref{claim1}, we have
\begin{align*}
  v_x&:= \max_{S\subseteq [n], |S|\le s} \sqrt{\bm a_{S}^{\top} (\bm B_{S, S})^{\dag} \bm a_{S}}= \sqrt{\bm a_{S^*}^{\top} (\bm B_{S^*, S^*})^{\dag} \bm a_{S^*}}.
\end{align*}
%where the first equation is from the combinatorial reformulation of problem \eqref{eq:sepx} by leveraging the result in the claim above. 

On the other hand,   the  Lagrangian dual of  problem  \eqref{eq:sepx} can be written as
    \begin{align*}
       v_x&\le \min_{\lambda \in \Re_+} \max_{\bm x\in \Re^n} \{\bm a^{\top} \bm x + \lambda - \lambda \bm x^{\top} \bm B \bm x:   \|\bm x\|_0 \le s_1  \}\\
       & =  \min_{\lambda \in \Re_+} \max_{S\subseteq [n], |S|\le s}  \lambda + \frac{\bm a_{S}^{\top} (\bm B_{S, S})^{\dag} \bm a_{S}}{4\lambda}\\
      & \le  \max_{S\subseteq [n], |S|\le s}  \lambda^* + \frac{\bm a_{S}^{\top} (\bm B_{S, S})^{\dag} \bm a_{S}}{4\lambda^*} = \sqrt{\bm a_{S^*}^{\top} (\bm B_{S^*, S^*})^{\dag} \bm a_{S^*}} \le v_x,
    \end{align*}
where the first equation is due to \Cref{claim1}, the second inequality is by plugging the feasible solution $\lambda^* = \frac{\sqrt{\bm a_{S^*}^{\top} (\bm B_{S^*, S^*})^{\dag} \bm a_{S^*}}}{2}$ into minimization problem, and the last equation is from the optimality of subset $S^*$. Since both left-hand and right-hand sides above equal $v_x$, the strong duality of problem \eqref{eq:sepx} holds, and all the inequalities above must attain the equalities. That is, problem \eqref{eq:sepx} is equivalent to
    \begin{align*}
       v_x&= \min_{\lambda \in \Re_+} \max_{\bm x\in \Re^n} \{\bm a^{\top} \bm x + \lambda - \lambda \bm x^{\top} \bm B \bm x:   \|\bm x\|_0 \le s_1  \}.
    \end{align*}
Since the outer minimization is a one-dimensional convex program that can be solved efficiently, as a result,  for any given $\lambda>0$, the inner maximization is equivalent to solving 
  \begin{align}\label{eq:sparse_scca_np}
    \max_{\bm x \in \Re^n}\{\bm a^{\top} \bm x  - \lambda\bm x^{\top} \bm B \bm x: \|\bm x\|_0 \le s_1 \}.
\end{align}

Next, let us consider the NP-hard sparse regression problem (see, e.g., \cite{natarajan1995sparse}), which admits
    \begin{align}\label{eq:sparse}
   \min_{\bm \beta\in \Re^n} \left\{\|\bm v -\bm U\bm x \|_2^2 : \|\bm x \|_0 \le s\right\} \Longleftrightarrow    \max_{\bm x\in \Re^n} \left\{2\bm v^{\top} \bm U\bm x- \bm x^{\top}\bm U^{\top} \bm U \bm \beta: \|\bm x\|_0 \le s\right\},
\end{align}
where data matrix $\bm U$ consists of observations of $n$ variables and vector $\bm v$ denotes the corresponding response variables.

Suppose that in the problem \eqref{eq:sparse_scca_np}, let us define $\lambda \bm B = \bm U^{\top} \bm U$ and $\bm a =2 \bm U^{\top}\bm v$. Then using the singular value decomposition of matrix $\bm U$, we see that the following equation still holds. $$\bm a_S - \bm B_{S,S}\bm B_{S,S}^{\dag}\bm a_S= \bm 0, \forall S\subseteq [n].$$
  Thus, for any given $\lambda>0$, the maximization problem \eqref{eq:sparse_scca_np} is equivalent to the sparse regression problem \eqref{eq:sparse}. This shows that problem \eqref{eq:sepx} is NP-hard.

   % Suppose that  $\lambda \bm B:=\bm Q \bm\Lambda\bm Q^{\top}$ denotes the eigendecomposition and the diagonal matrix $\bm \Lambda$ consists of all positive eigenvalues. 
   % By setting  $\bm U:=\sqrt{\bm\Lambda}\bm Q^{\top}$, we have that $\lambda \bm B = \bm U^{\top} \bm U$. In addition,  we construct the vector $\bm v$ based on the equation $2 \bm U^{\top}\bm v=\bm a$, implying that
   % $$2\bm Q \sqrt{\bm\Lambda} \bm v =  \bm a     \Longrightarrow \bm v = \frac{\sqrt{\bm \Lambda^{-1}} \bm Q^{\top} \bm a}{2},$$
   % where the result is from the fact $\bm Q^{\top} \bm Q = \bm I$.
    
Similarly, the maximization problem over $\bm y$ in \eqref{eq:sep} can also be reduced to the sparse regression problem.
    \qed

% \section*{Appendix B: Supplementary computational results}

% \subsection{Real data}
% The real data consist of samples of two variables: $2149\times 89$ and $19672\times 89$, where $n=2149$ and $m=19672$. The rank is at most $89$.

% \begin{table}[htb]
% 		\centering
%     \scriptsize
%   \caption{Real Data}
% 		\setlength{\tabcolsep}{3pt}\renewcommand{\arraystretch}{1.3}
% 		\begin{tabular}{c c c c c|r r| r r| r  r | r r}
% 			\hline 
% 		 & & & & & \multicolumn{2}{c|}{Greedy} & \multicolumn{2}{c|}{Local Search} & \multicolumn{2}{c|}{SDP relaxtion} & \multicolumn{2}{c}{Branch-and-cut}
% 			\\ \cline{6-13}
% 		$n$ & $m$ & $s_1$ &	 $s_2$& $r$ & \multicolumn{1}{c}{LB}  & \multicolumn{1}{c|}{time(s)}   &  \multicolumn{1}{c}{LB}  & \multicolumn{1}{c|}{time(s)}  &\multicolumn{1}{c}{UB} &\multicolumn{1}{c|}{time(s)} &\multicolumn{1}{c}{$v^*$} &\multicolumn{1}{c}{time(s)}   \\
% 			\hline
% 50 & 50 & 10 & 10 & 50 & 0.9712 & 1 & 0.9712 & 1 & 0.9999 & 12\\
% 50 & 100 & 10 & 10 & 50 & 0.9784 & 1 & 0.9791 & 1 & 1.0000 & 63\\
% 100 & 100 & 10 & 15 & 89 & 0.9882 & 1 & 0.9882 & 1 &  1.0000 & 1077\\
% 100 & 150 & 10 & 20 & 89& 0.9919 & 1 & 0.9924 & 3 \\
% 150 & 150 & 20 & 20 & 89&0.9967 & 1 & 0.9980 & 28 \\
%    \hline
% 		\end{tabular}%
%   \label{table2}
% \end{table}

\section*{Appendix B: Implementations of greedy and local search algorithms}
This section presents the detailed implementations of greedy and local search algorithms based on the formulation \eqref{scca:com}.
\begin{algorithm}
	\caption{Greedy algorithm for SCCA \eqref{scca:com}}
	\label{algo:svd_greedy}
	\begin{algorithmic}[1]
		\State \textbf{Input:} Matrices $\bm A \in \Re^{m \times n}$, $\bm B \in \S_+^m$, $\bm C \in \S_+^m$ and  integers   $s_1\in [n]$, $s_2\in [m]$
		
		\State Compute $(i^*, j^*) \in \argmax_{i \in[m], j \in[n] }\sqrt{(B_{ii})^{\dag}}  A_{ij} \sqrt{(C_{jj})^{\dag}} $
		\State  Define subsets $\hat{S}_1 := \{i^*\}$ and $\hat{S}_2 := \{j^*\}$

		\For{$\ell = 2, \cdots, \max\{s_1,s_2\}$}
		\If{$\ell \le \min\{s_1,s_2\}$}
		
		\State 
		$i^* \in \argmax_{i \in[n]\setminus \hat{S}_1  }
	\sigma_{\max}\left(\sqrt{(\bm B_{\hat{S}_1 \cup\{i\}, \hat{S}_1 \cup\{i\}})^{\dag}} \bm A_{\hat{S}_1 \cup\{i\}, \hat{S}_2} \sqrt{(\bm C_{\hat{S}_2, \hat{S}_2})^{\dag}}\right)
		$
		\State Update  $\hat{S}_1 := \hat{S}_1\cup \{i^*\}$ 
		
		\State $j^* \in \argmax_{j\in[m]\setminus \hat{S}_2  }
	\sigma_{\max}\left(\sqrt{(\bm B_{\hat{S}_1 , \hat{S}_1})^{\dag}} \bm A_{\hat{S}_1, \hat{S}_2 \cup\{j\}} \sqrt{(\bm C_{\hat{S}_2\cup\{j\}, \hat{S}_2\cup\{j\}})^{\dag}}\right)$
		
		\ElsIf{$s_1\leq s_2$}
		\State 
		$j^* \in \argmax_{j\in[m]\setminus \hat{S}_2  }
	\sigma_{\max}\left(\sqrt{(\bm B_{\hat{S}_1 , \hat{S}_1})^{\dag}} \bm A_{\hat{S}_1, \hat{S}_2 \cup\{j\}} \sqrt{(\bm C_{\hat{S}_2\cup\{j\}, \hat{S}_2\cup\{j\}})^{\dag}}\right)$ 
		\State  Update $\hat{S}_2:= \hat{S}_2 \cup \{ j^*\}$
		\Else
		\State 
		$i^* \in \argmax_{i \in[n]\setminus \hat{S}_1  }
	\sigma_{\max}\left(\sqrt{(\bm B_{\hat{S}_1 \cup\{i\}, \hat{S}_1 \cup\{i\}})^{\dag}} \bm A_{\hat{S}_1 \cup\{i\}, \hat{S}_2} \sqrt{(\bm C_{\hat{S}_2, \hat{S}_2})^{\dag}}\right)
		$
		\State  Update  $\hat{S}_1 := \hat{S}_1\cup \{i^*\}$
		\EndIf
		\EndFor
		
		\State  \textbf{Output:} $\hat{S}_1 , \hat{S}_2 $%$(\bar{x},\bar{w})$ 
	\end{algorithmic}
\end{algorithm}

\begin{algorithm}[htb]
	\caption{Local search algorithm for SSVD \eqref{scca:com}}
	\label{algo:svd_localsearch}
	\begin{algorithmic}[1]
		\State \textbf{Input:} Matrices $\bm A \in \Re^{m \times n}$, $\bm B \in \S_+^m$, $\bm C \in \S_+^m$ and  integers   $s_1\in [n]$, $s_2\in [m]$
		%		\State Let $\bm{A}=\bm{C}^{\top}\bm{C}$ denote its Cholesky factorization where $\bm{C} \in \Re^{d\times n}$
		%		\State Let $\bm{c}_i \in \Re^d$ denote the $i$-th column vector of matrix $\bm{C}$ for each $i \in [n]$
		\State Initialize $(\hat{S}_1 ,\hat{S}_2)$ as the output of greedy Algorithm~\ref{algo:svd_greedy}
		\State \textbf{do}
		\For {each pair {$(i_1,j_1) \in \hat{S}_1 \times ([n]\setminus \hat{S}_1)$}}
		\If{$\sigma_{\max}\left(\sqrt{(\bm B_{\hat{S}_1\cup \{j_1\} \setminus \{i_1\}, \hat{S}_1\cup \{j_1\} \setminus \{i_1\}})^{\dag}} \bm A_{\hat{S}_1\cup \{j_1\} \setminus \{i_1\}, \hat{S}_2} \sqrt{(\bm C_{\hat{S}_2, \hat{S}_2})^{\dag}}\right)> \sigma_{\max}\left(\sqrt{(\bm B_{\hat{S}_1, \hat{S}_1})^{\dag}} \bm A_{\hat{S}_1, \hat{S}_2} \sqrt{(\bm C_{\hat{S}_2, \hat{S}_2})^{\dag}}\right)$}
		\State Update $\hat{S}_1 := \hat{S}_1 \cup \{j_1\} \setminus \{i_1\}$
		\EndIf
		\EndFor
		
		\For{each pair {$(i_2,j_2) \in \hat{S}_2 \times ([m]\setminus\hat{S}_2)$}}
		\If{$\sigma_{\max}\left(\sqrt{(\bm B_{\hat{S}_1\cup \{j_1\} \setminus \{i_1\}, \hat{S}_1\cup \{j_1\} \setminus \{i_1\}})^{\dag}} \bm A_{\hat{S}_1\cup \{j_1\} \setminus \{i_1\}, \hat{S}_2} \sqrt{(\bm C_{\hat{S}_2, \hat{S}_2})^{\dag}}\right)> \sigma_{\max}\left(\sqrt{(\bm B_{\hat{S}_1, \hat{S}_1})^{\dag}} \bm A_{\hat{S}_1, \hat{S}_2} \sqrt{(\bm C_{\hat{S}_2, \hat{S}_2})^{\dag}}\right)$}
		\State Update  $\hat{S}_2 := \hat{S}_2 \cup \{j_2\} \setminus \{i_2\}$
		\EndIf
		\EndFor
		\State\textbf{while} {there is still an improvement}
		\State \textbf{Output:} $\hat{S}_1$, $\hat{S}_2$%$(\bar{x},\bar{w})$
	\end{algorithmic}
\end{algorithm}

\end{document}